\documentclass[a4paper, 10pt]{article} \usepackage{amsmath}

\usepackage{amscd}

\usepackage{amssymb}

\usepackage{amsthm}    \newtheorem{thm}{Theorem}[section]

\newtheorem{prop}[thm]{Proposition}

\newtheorem{lem}[thm]{Lemma}

\newtheorem{cor}[thm]{Corollary}  \theoremstyle{definition}

\newtheorem{df}[thm]{Definition}  \theoremstyle{definition}

\newtheorem{rem}{Remark}        \theoremstyle{plain}

 \theoremstyle{definition}

\newtheorem{ex}[thm]{Example}  \def\CC{\Bbb{C}}

\def\RR{\Bbb{R}} 

\def\CCI{\overline{\CC}}    \def\NN{\Bbb{N}}

\def\B1{{\rm\kern.32em\vrule  width.12em    height1.4ex

depth-.05ex\kern-.28em 1}}

\begin{document}

\title{Dimensions of Julia sets of\\ expanding rational
semigroups\footnote{
Published in Kodai Mathematical Journal, 
Vol. 28, No.2 (2005), pp. 390-422.  
 The title was changed. 
The first version (of which title was "The dimensions of Julia sets of 
expanding rational semigroups") was written on 
May 27,2004.
Keywords: expanding rational semigroup,
Hausdorff dimension of a Julia set,
thermodynamic formalism. 
}}
\author{Hiroki Sumi\
\\ Department of Mathematics, Graduate School of Science,\\ 
Osaka University,\\ 
1-1,\ Machikaneyama,\ Toyonaka,\ Osaka,\ 560-0043,\   
Japan\\ E-mail:    sumi@math.sci.osaka-u.ac.jp\\
http://www.math.sci.osaka-u.ac.jp/$\sim $sumi/welcomeou-e.html\\  
 \\ 
{\em Dedicated to the memory of Professor Nobuyuki Suita}}

\maketitle

\begin{abstract}

We estimate the upper box
and Hausdorff dimensions of the Julia set of
an expanding semigroup generated
by finitely many rational functions,
using the thermodynamic formalism in ergodic theory.
Furthermore, we show Bowen's formula, and the existence and uniqueness
of a conformal measure, for a finitely generated expanding
semigroup satisfying the open set condition.
\end{abstract}
\section{Introduction}
For a Riemann surface $S$, let End($S$) denote the set
of all holomorphic endomorphisms of $S$. This is a semigroup
whose semigroup operation is the functional composition.
A {\bf rational semigroup} is a subsemigroup of End($\CCI $)
without any constant elements.  We say that a rational
semigroup $G$ is a {\bf polynomial semigroup} if each
element of $G$ is a polynomial. Research on the dynamics of
rational semigroups was initiated by
A. Hinkkanen and G.J. Martin (\cite{HM1}),
who were interested in the role of the
dynamics of polynomial semigroups while studying
various one-complex-dimensional
moduli spaces for discrete groups,
and
F. Ren's group(\cite{ZR}, \cite{GR}).
For references on research into rational semigroups,
see \cite{HM1}, \cite{HM2}, \cite{HM3}, \cite{ZR}, \cite{GR},
\cite{SSS}, 
\cite{Bo}, \cite{St1}, \cite{St2}, \cite{St3},
\cite{S1}, \cite{S2}, \cite{S3}, \cite{S4}, \cite{S5},
\cite{S6}, and \cite{S7}.
The research on the dynamics of rational semigroups
can be considered a generalization of
studies of both the {\bf iteration of rational functions} and
{\bf self-similar sets constructed using iterated function systems of
some similarity transformations in $\RR ^{2}$} in fractal geometry.
In both fields, the estimate of
the upper(resp.lower) box dimension,
which is denoted by $\overline{\dim }_{B}$
(resp.$\underline{\dim } _{B}$), and the
Hausdorff dimension, which is denoted by
$\dim _{H}$, of the invariant sets
(Julia sets or attractors) has been of
great interest and has been investigated for a long time.
In this paper, we consider the following: 
 For a rational semigroup $G$, 
 We set
\[ F(G) = \{x\in \CCI \mid G \mbox{ is normal in a neighborhood of $x$}\},
\ J(G) = \CCI \setminus F(G) .\] \(F(G)\) is called the
{\bf Fatou set} for $G$ and \(J(G)\) is called the {\bf
Julia set} for $G$. We use $\langle f_{1}, f_{2}, \cdots \rangle$ to denote the
rational semigroup generated by the family $\{f_{i}\}$.
For a finitely generated rational semigroup 
$G=\langle f_{1},\cdots ,f_{m}\rangle $, 
we set $\Sigma _{m}=\{1,\cdots ,m\}^{\NN}$ (this is a 
compact metric space) 
and we use $\sigma :\Sigma _{m}\rightarrow \Sigma _{m}$ to denote the shift map,
which is
$(w_{1},w_{2},\cdots ) \mapsto (w_{2},w_{3},\cdots )$ for
$w=(w_{1},w_{2},w_{3},\cdots)\in \Sigma _{m}$.
We define the map
$\tilde{f} :\Sigma _{m}\times \CCI \rightarrow \Sigma _{m}\times \CCI$ using
\[\tilde{f}((w,x))=(\sigma w, f_{w_{1}}(x)).\] 
We call map $\tilde{f}$ the {\bf skew product map associated with
the generator system $\{ f_{1},\cdots ,f_{m}\}$.}
 For each $w\in \Sigma _{m}$, we use $F_{w}$ to denote
the set of all the points $x\in \CCI$ that satisfy the fact that there exists
an
open neighborhood $U$ of $x$ such that the family 
$\{ f_{w_{n}}\circ \cdots \circ f_{w_{1}}\} _{n}$ is normal in $U$. 
We set $J_{w}=\CCI \setminus F_{w}$ and
$\tilde{J}_{w}=\{w\} \times J_{w}$. Moreover, we set
$$ \tilde{J}(\tilde{f})=\overline{\bigcup _{w\in \Sigma _{m}}\tilde{J}_{w}},\ 
\tilde{F}(\tilde{f})=(\Sigma _{m}\times \CCI)\setminus\tilde{J}(\tilde{f}),$$
where the closure is taken in the product space 
$\Sigma _{m}\times \CCI $ (this is a compact metric space).
We call
$\tilde{F}(\tilde{f})$ the {\bf Fatou set} for $\tilde{f}$ and
$\tilde{J}(\tilde{f})$ the {\bf Julia set} for $\tilde{f}$.
For each $(w, x) \in \Sigma _{m}\times \CCI$ and $n \in \NN$ we set
\[(\tilde{f}^{n})'((w,x))= (f_{w_{n}}\cdots f_{w_{1}})'(x).\]
Furthermore, we denote the first (resp. second) projection by
$\pi :\Sigma _{m}\times \CCI \rightarrow \Sigma _{m}$
(resp. $\pi _{\CCI}:\Sigma _{m}\times \CCI \rightarrow
\CCI$).
 We say that a finitely generated 
rational semigroup $G=\langle f_{1},\cdots ,f_{m}\rangle $ is 
an {\bf expanding rational semigroup}
if $J(G)\neq \emptyset$ and
the skew product map $\tilde{f}: \Sigma _{m}\times \CCI \rightarrow
\Sigma _{m}\times \CCI$ associated with the generator system
$\{f_{1},\cdots ,f_{m}\}$ is expanding along fibers, i.e.,
there exists a positive constant $C$ and a constant $\lambda >1$
such that for each $n\in \NN$,
$$\inf\limits _{z\in \tilde{J}(\tilde{f})} \| (\tilde{f}^{n})'(z)\| \geq
C\lambda ^{n},$$
where we use $\| \cdot \|$ to denote the norm of the derivative with respect
to the spherical metric.

For a general rational semigroup $G$ and  
a non-negative number $t$, we say
that a Borel probability measure $\tau$ on $\CCI$ is {\bf $t$-subconformal} (for $G$)
 if for
each $g\in G$ and for each Borel measurable set $A$ in $\CCI$,
$\tau (g(A))\leq \int _{A} \| g'\| ^{t}d\tau$.
Moreover, we set 
$s(G)=\inf \{t \mid \exists \tau : t \mbox{-subconformal measure}\} 
.$

 Furthermore,\ we say that a Borel probability measure
$\tau$ on $J(G)$ is
$t$-{\bf conformal} (for $G$)
if
for any Borel set $A$ and
$g\in G,$ if $A,g(A)\subset J(G)$ and
$g:A\rightarrow g(A)$ is injective, then
$\tau (g(A))=\int _{A}\| g'\| ^{t}\ d\tau .$

 For any $s\geq 0$ and $x\in \CCI$,
we set $S(s,x)=\sum _{g\in G}\sum _{g(y)=x}$
$\|g'(y)\| ^{-s}$. Moreover, we set
$S(x)=\inf \{s\geq 0\mid S(s,x)<\infty \}$
(If no $s$ exists with $S(s,x)<\infty$,
then we set $S(x)=\infty$). We set 
$s_{0}(G)=\inf \{ S(x)\mid x\in \CCI \} .$ 
Note that if $G$ has only countably many elements, 
then $s(G)\leq s_{0}(G)$ (Theorem 4.2 in \cite{S2}.)

 Then, under the above notations,\ we show the following:
\begin{thm}({\bf Main Theorem A})
\label{mainthA}
Let
$G=\langle f_{1},\cdots ,f_{m}\rangle $ be
a finitely generated expanding rational semigroup.
Let $\tilde{f}:\Sigma _{m}\times \CCI 
\rightarrow \Sigma _{m}\times \CCI $ be the skew product map 
associated with
$\{ f_{1},\cdots ,f_{m}\} $.
Then, there exists a unique zero
$\delta $ of the function:
$P(t):=P(\tilde{f}|_{\tilde{J}(\tilde{f})}, t\tilde{\varphi })$,
where $\tilde{\varphi }$ is the function on
$\tilde{J}(\tilde{f})$ defined by:
$\tilde{\varphi }((w,x))=-\log (\| (f_{w_{1}})'(x)\| )$ 
for $(w,x)=((w_{1},w_{2},\cdots ),x)\in \tilde{J}(\tilde{f})$ and
$P(\ ,\ )$ denotes the pressure.
Furthermore, $\delta $
satisfies the fact that there exists a unique probability
measure $\tilde{\nu }$ on $\tilde{J}(\tilde{f}) $
such that $M_{\delta }^{\ast }\tilde{\nu }=\tilde{\nu } $,
where $M_{\delta } $ is an operator on
$C(\tilde{J}(\tilde{f}))$ (the space of continuous functions 
on $\tilde{J}(f)$) defined by
$$ M_{\delta }\psi ((w,x))=\sum _{\tilde{f} ((w',y))=(w,x)}
\frac{\psi ((w',y))}{\| (f_{w_{1}'})'(y)\| ^{\delta }},$$
where $w'=(w_{1}',w_{2}',\cdots )\in \Sigma _{m}.$ 
Moreover, $\delta $ satisfies
$$ \overline{\dim }_{B}(J(G))\leq s(G)\leq s_{0}(G)
\leq \delta =\frac{h_{\alpha \tilde{\nu }}(\tilde{f})}
{-\int _{\tilde{J}(\tilde{f})}\tilde{\varphi }\alpha d\tilde{\nu }}\leq 
\frac{\log (\sum _{j=1}^{m}\deg (f_{j}))}
{-\int _{\tilde{J}(\tilde{f})}\tilde{\varphi} \alpha d\tilde{\nu}},$$
where $\alpha =\lim _{l\rightarrow \infty}M_{\delta}^{l}(1)$ 
and we denote the metric entropy of
$(\tilde{f},\alpha \tilde{\nu})$ by
$h_{\alpha \tilde{\nu}}(\tilde{f})$.
The support for $\nu :=(\pi _{\CCI})_{\ast}\tilde{\nu}$ equals $J(G)$.

 Furthermore, let
$A(G)=\overline{\cup _{g\in G}
g(\{ x\in \CCI \mid \exists h\in G, h(x)=x, |h'(x)|<1\} )} $ and 
$P(G)=\overline{\bigcup _{g\in G}\{ \mbox{\em all critical values of } g\} }.$ 
Then, $A(G)\cup P(G)\subset F(G)$ and for each
$x\in \CCI \setminus (A(G)\cup P(G))$, we have
$\delta$ is equal to: 
$$\inf \{t\geq 0\mid 
\sum _{n\in \NN}\sum _{(w_{1},\cdots ,w_{n})\in 
\{1,\cdots ,m\} ^{n}}\sum 
_{(f_{w_{1}}\cdots 
f_{w_{n}})(y)=x} 
\| (f_{w_{1}}\cdots
f_{w_{n}})'(y)\| ^{-t}<\infty \} .$$
\end{thm}
\begin{thm}({\bf Main Theorem B})
\label{mainthB}
Let
$G=\langle f_{1},\cdots ,f_{m}\rangle$ be
a finitely generated expanding rational semigroup.
Suppose that there exists a non-empty
open set U in $\CCI$ such that
$f_{j}^{-1}(U)\subset U$ for each
$j=1,\cdots ,m$ and
$\{f_{j}^{-1}(U)\} _{j=1}^{m}$ are mutually
disjoint. Then, we have the following:
\begin{enumerate}
\item
$\dim _{H}(J(G))=\overline{\dim}_{B}(J(G))=
s(G)=s_{0}(G)=\delta$, where $\delta$ denotes the number
in Theorem~\ref{mainthA}.
\item
$\nu:=(\pi _{\CCI})_{\ast}\tilde{\nu}$ is the
unique $\delta$-conformal measure, where $\tilde{\nu}$
is the measure in Theorem~\ref{mainthA}.
Furthermore, $\nu$ satisfies the fact that there exists a positive constant
$C$ such that for any $x\in J(G)$ and any positive number
$r$ with $r<$ {\em diam} $\CCI$, we have
$$C^{-1}\leq \frac{\nu (B(x,r))}{r^{\delta}}\leq C.$$
\item $\nu$ satisfies $\nu (f_{i}^{-1}(J(G))\cap
f_{j}^{-1}(J(G)))=0$, for each
$i,j\in \{1,\cdots ,m\}$ with $i\neq j$.
Furthermore, for each $(i,j)$ with $i\neq j$,
we have
$f_{i}^{-1}(J(G))\cap f_{j}^{-1}(J(G))$ is
nowhere dense in $f_{j}^{-1}(J(G))$.
\item $0<H^{\delta}(J(G))<\infty $, where
$H^{\delta}$ denotes the $\delta$-dimensional Hausdorff measure
with respect to the spherical metric. Furthermore,
we have $\nu =\frac{H^{\delta}|_{J(G)}}{H^{\delta}(J(G))}$.
\item If there exists a $t$-conformal measure $\tau $,
then $t=\delta$ and $\tau =\nu$.
\item For any $x\in \CCI \setminus (A(G)\cup P(G))$,
we have
$$\dim _{H}(J(G))=\delta =
\inf \{t\geq 0\mid \sum _{g\in G}\sum _{g(y)=x}
\| g'(y)\| ^{-t}<\infty \} .$$
\end{enumerate}
\end{thm}
\begin{rem}
In \cite{S6}, it is shown that
if $G=\langle f_{1}\cdots ,f_{m}\rangle$
is expanding and there exists a non-empty open set
$U$ such that $f_{j}^{-1}(U)\subset U$ for each
$j=1,\cdots ,m$, $\{f_{j}(U)\} _{j}$ are
mutually disjoint and $\overline{U}\neq J(G)$, then
$J(G)$ is porous and $\overline{\dim}_{B}(J(G))<2$.
\end{rem}
\begin{rem}
In addition to the assumption
of Main Theorem B,
if $J(G)\subset \CC$, then we can also
show a similar result for the Euclidean metric.
\end{rem}
For the precise notation,
see the following sections.
The proof of Main Theorem A is given in
section~\ref{Hausdimex} and the proof
of Main Theorem B is given in
section~\ref{secmainB}.
The existence of a subconformal or conformal measure
is deduced by
applying some of the results in \cite{W1} and
the thermodynamic formalism in ergodic theory
to the skew product map associated with the generator system.
Since generator maps are not injective in general and
we do not assume the ``cone condition'' 
(the existence of uniform cones) for the boundary of
the open set,
much effort is needed to estimate
$\nu (B(x,r))$ in Main Theorem B. Indeed,   
we cut the closure of the open set 
into small pieces $\{ K_{j}\} $, and 
for a fixed $s\in \NN $, 
let ${\cal K}$ be the set of all 
$(\gamma , k_{j})$ that satisfies that 
$\gamma $ is a well defined inverse branch 
of $(f_{w_{1}}\circ \cdots f_{w_{u}})^{-1}$ defined 
on $K_{j}$  
for some $(w_{1},\cdots ,w_{u})\in \{ 1,\cdots ,m\} ^{u}$ 
with $u\leq s.$ 
Then we introduce an equivalence class $``\sim "$ 
in a subset $\Gamma $ of ${\cal K}$, and 
an order $``\leq "$ in $\Gamma /\! \sim .$ 
We obtain an upper estimate of 
the cardinality of the set of all minimal elements 
of $(\Gamma /\! \sim ,\ \leq )$ by a constant independent 
of $r$ and $x$, 
which gives us the key to estimate $\nu (B(x,r)).$ 

 Note that in \cite{MU1}, it was discussed the case 
 in which there are infinitely many injective generator maps 
 and the boundary of the open set satisfies the cone condition. 

 The uniqueness of a conformal measure $\tau$ is deduced from
some results in \cite{W1} and an estimate of
$\tau (B(x,r))$. Note that our definition of conformal measure
differs from that of \cite{MU1} and \cite{MU2}.
In this paper, we do not require the separating condition
for the definition of conformal measure.

\section{Preliminaries}
In this section, we give the notation and definitions
for rational semigroups and the associated skew products
that we need to give our main result.
\subsection{Rational semigroups}
\label{rs}
We use the definition in \cite{S5}.
\begin{df}
Let $G$ be a rational semigroup. We set
\[ F(G) = \{z\in \CCI \mid G \mbox{ is normal in a neighborhood of $z$}\},
\ J(G) = \CCI \setminus F(G) .\] \(F(G)\) is called the
{\bf Fatou set} for $G$ and \(J(G)\) is called the {\bf
Julia set} for $G$. The backward orbit $G^{-1}(z)$ of $z$ and
the
{\bf set of exceptional points} $E(G)$ are defined
by:
$ G^{-1}(z) =\cup _{g\in G} g^{-1}(z)$ and
$ E(G)= \{z\in \CCI \mid \sharp G^{-1}(z) \leq2 \}$.
For any subset $A$ of $\CCI$, we set
$G^{-1}(A)=\cup _{g\in G}g^{-1}(A)$.
We use $\langle f_{1}, f_{2}, \cdots \rangle$ to denote the
rational semigroup generated by the family $\{f_{i}\}$.
For a rational map $g$, we use $J(g)$ to denote the Julia set
of dynamics of $g$.
\end{df}
For a rational semigroup $G$, for each $f\in G$, we have
$f(F(G))\subset F(G)$ and $f^{-1}(J(G))\subset J(G)$. Note that we do not
have this equality hold in general.
If  $\sharp J(G)\geq 3 $, then \( J(G) \) is a
perfect set,
$\sharp E(G) \leq 2$,
\(J(G)\) is the smallest closed backward invariant set containing at least
three points, and $J(G)$ is the closure of the union of all
repelling fixed points of elements of $G$,
which implies that
$J(G)=\overline{\bigcup _{g\in G}
J(g)}$.
If a point \(z\) is not in \(E(G),\) then for every \(x\in J(G)\),
$x\in \overline{G^{-1}(z)}.$ In particular,
if \ \(z\in J(G)\setminus E(G)\), then
$\overline{G^{-1}(z)}=J(G)$.
Further, for a finitely generated
rational semigroup $G=\langle
f_{1}, \cdots , f_{m}\rangle $, if
we use $G_{n}$ to denote the subsemigroup of
$G$ that is generated by $n$-products
of generators $\{f_{j} \}$, then
$J(G_{n})=J(G)$.
For more precise statements, see
Lemma 2.3 in \cite{S5}, for which the proofs are
based on \cite{HM1} and \cite{GR}.
Furthermore, if $G$ is generated by a precompact
subset $\Lambda$ of End$(\CCI)$, then
$J(G)=\overline{\bigcup _{f\in \Lambda}f^{-1}(J(G))}=
\bigcup _{h\in \overline{\Lambda}}h^{-1}(J(G))$.
In particular, if $\Lambda$ is compact, then we have
$J(G)= \bigcup _{f\in \Lambda}f^{-1}(J(G))$(\cite{S3}).
We call this property of a Julia set
the {\bf backward self-similarity}.
\begin{rem}
Using the backward self-similarity, research
on the Julia sets of rational semigroups
may be considered a generalization
of research on self-similar sets constructed
using some similarity transformations from $\CC$ to itself,
which can be regarded as the Julia sets of some
rational semigroups. It is easily seen that
the Sierpi\'{n}ski gasket is the Julia set
of a rational semigroup
$G=\langle f_{1},f_{2},f_{3}\rangle$ where
$f_{i}(z)=2(z-p_{i})+p_{i},i=1,2,3$ with
$p_{1}p_{2}p_{3}$ being a regular triangle.
\end{rem}
\subsection{Associated skew products}
We use the notation in \cite{S5}.
Let $m$ be a positive integer. We use $\Sigma _{m}$ to denote the one-sided wordspace that is
$\Sigma _{m}=\{1,\cdots ,m\}^{\NN}$
and use $\sigma :\Sigma _{m}\rightarrow \Sigma _{m}$ to denote the shift map,
which is
$(w_{1},\cdots)\mapsto (w_{2},\cdots)$ for
$w=(w_{1},w_{2},w_{3},\cdots)\in \Sigma _{m}$.
For any $w,w'\in \Sigma _{m}$,
we set $d(w,w'):=\sum _{n=1}^{\infty}(1/2^{n})\cdot
c(w_{k},w'_{k})$, where
$c(w_{k},w'_{k})=0$ if
$w_{k}=w'_{k}$ and $c(w_{k},w'_{k})=1$ if $w_{k}\neq w'_{k}$.
Then, $(\Sigma _{m},d)$ is a compact metric space.
Furthermore, the dynamics of 
$\sigma :\Sigma _{m}\rightarrow \Sigma _{m}$ is 
expanding with respect to this metric $d$. That is,
each inverse branch $\sigma _{j}^{-1}$ of $\sigma ^{-1}$
on $\Sigma _{m}$,
which is defined by $\sigma _{j}^{-1}((w_{1}, w_{2},\cdots))
=(j,w_{1},w_{2},\cdots)$ for $j=1,\cdots ,m$,
satisfies
$d(\sigma _{j}^{-1}(w), \sigma _{j}^{-1}(w'))\leq
(1/2)\cdot d(w,w')$.

Let $G= \langle f_{1},f_{2},\cdots ,f_{m} \rangle$ be a finitely generated rational
semigroup. We define the map
$\tilde{f} :\Sigma _{m}\times \CCI \rightarrow \Sigma _{m}\times \CCI$ using
\[\tilde{f}((w,x))=(\sigma w, f_{w_{1}}(x)).\] 
We call map $\tilde{f}$ the {\bf skew product map associated with
the generator system $\{ f_{1},\cdots ,f_{m}\}$.}
$\tilde{f}$ is finite-to-one and an open map. We hold that point $(w,x)\in
\Sigma _{m}\times \CCI$ satisfies $f_{w_{1}}'(x)\neq 0$ if and only if
$\tilde{f}$ is a homeomorphism in a small neighborhood of $(w,x)$. Hence, the map
$\tilde{f}$ has infinitely many critical points in general.
\begin{df}
For each $w\in \Sigma _{m}$, we use $F_{w}$ to denote
the set of all the points $x\in \CCI$ that satisfy the fact that there exists
an
open neighborhood $U$ of $x$ such that the family 
$\{f_{w_{n}}\circ \cdots \circ f_{w_{1}}\} _{n}$ is normal in $U$. 
We set $J_{w}=\CCI \setminus F_{w}$ and
$\tilde{J}_{w}=\{w\} \times J_{w}$. Moreover, we set
$$ \tilde{J}(\tilde{f})=\overline{\bigcup _{w\in \Sigma _{m}}\tilde{J}_{w}},
\  
\tilde{F}(\tilde{f})=(\Sigma _{m}\times \CCI)\setminus\tilde{J}(\tilde{f}),$$
where the closure is taken in the product space 
$\Sigma _{m}\times \CCI .$
We often write $\tilde{F}(\tilde{f})$ as
$\tilde{F}$ and $\tilde{J}(\tilde{f})$ as $\tilde{J}$. We call
$\tilde{F}(\tilde{f})$ the {\bf Fatou set} for $\tilde{f}$ and
$\tilde{J}(\tilde{f})$ the {\bf Julia set} for $\tilde{f}$. Here, we remark
that $\bigcup _{w\in \Sigma _{m}}\tilde{J}_{w}$ may not be compact
in general. That is why we consider the closure of that set in
$\Sigma _{m}\times \CCI$ (this is a compact space) concerning
the definition of the Julia set for $\tilde{f}$.

For each $(w, x) \in \Sigma _{m}\times \CCI$ and $n \in \NN$ we set
\[(\tilde{f}^{n})'((w,x))= (f_{w_{n}}\cdots f_{w_{1}})'(x).\]
Furthermore, we denote the first (resp. second) projection by
$\pi :\Sigma _{m}\times \CCI \rightarrow \Sigma _{m}$
(resp. $\pi _{\CCI}:\Sigma _{m}\times \CCI \rightarrow
\CCI$).
Note that we have
$\tilde{f}(\tilde{F}(\tilde{f}))=\tilde{f}^{-1}(\tilde{F}(\tilde{f}))=
\tilde{F}(\tilde{f}),\ \tilde{f}(\tilde{J}(\tilde{f}))=
\tilde{f}^{-1}(\tilde{J}(\tilde{f}))=\tilde{J}(\tilde{f})$ and
$\pi _{\CCI}(\tilde{J}(\tilde{f}))=J(G)$.
(For the fundamental properties of these sets,
see Proposition 3.2 in \cite{S5}. In addition, see \cite{S3}.)
\end{df}

\begin{df}
Let $G=\langle f_{1},\cdots ,f_{m}\rangle$ be
a finitely generated rational semigroup.
Let us fix the generator system
$\{f_{1},\cdots ,f_{m}\}$.
We set $f_{w}:=f_{w_{1}}\circ \cdots \circ f_{w_{k}}$
for any $w=(w_{1},\cdots ,w_{k})\in
\{1,\cdots ,m\} ^{k}$.
We set ${\cal W}=\cup _{n\in \NN}
\{1,\cdots ,m\} ^{n}\bigcup \Sigma _{m}$ and
set ${\cal W}^{\ast}=\cup _{n\in \NN}
\{1,\cdots ,m\} ^{n}$.
For any $w=(w_{1},w_{2},\cdots)\in {\cal W}$, we set
$|w|=n$ if $w\in \{1,\cdots ,m\} ^{n}$ and
$|w|=\infty$ if $w\in \Sigma _{m}$.
Furthermore, we set
$w|k:= (w_{1},\cdots ,w_{k})$, for any
$k\in \NN$ with $k\leq |w|$.
Moreover,
for any $w\in {\cal W}^{\ast}$,
we set $\Sigma _{m}(w):=\{w'\in \Sigma _{m}
\mid w'_{j}=w_{j}, j=1,\cdots ,|w|\} $.
For any $w^{1}\in {\cal W}^{\ast}$ and
$w^{2}\in {\cal W}$,
we set $w^{1}w^{2}=
(w_{1}^{1},\cdots ,w_{|w^{1}|}^{1},w_{1}^{2},w_{2}^{2},\cdots)
\in {\cal W}$.
\end{df}
\noindent {\bf Notation:}
Let $(X,d)$ be a metric space.
For any subset $A$ of $X$,
we set diam $A:=\sup \{d(x,y)\mid x,y\in A\}$.
Let $\mu$ be a Borel measure on $X$.
We use
supp $\mu$ to denote the support of $\mu$.
For any Borel set $A$ in $X$, we use
$\mu |_{A}$ to denote the measure on $A$
such that $\mu |_{A}(B)=\mu (B)$ for each
Borel subset $B$ of $A$.
We set $L^{1}(\mu)=\{\varphi :X\rightarrow \RR \mid
\int _{X}|\varphi |d \mu <\infty \}$,
with $L^{1}$ norm.
For any $\varphi \in L^{1}(\mu)$,
we sometimes use $\mu (\varphi)$ to mean
$\int _{X}\varphi \ d\mu$.
For
any $\varphi \in L^{1}(\mu)$,
we use $\varphi \mu$ to denote the measure
such that $(\varphi \mu)(A)=\int _{A}\varphi \ d\mu$
for any Borel set $A$.
We set $C(X)=\{\varphi :X\rightarrow \RR \mid
\mbox{continuous}\}$. (If $X$ is compact,
then $C(X)$ is the Banach space with the supremum norm.)
For any subset $A$ of $X$ and any $r>0$,
we set $B(A,r)=\{y\in X \mid d(y,A)<r\}$.
For any subset $A$ of $X$, we use
int $A$ to denote the interior of $A$.
\begin{rem}
\label{sande}
In this paper, we always use the
spherical metric on $\CCI$. However,
we note that conjugating a
rational semigroup $G$ by a M\"{o}bius
transformation, we may assume
that $J(G)\subset \CC$, and then
for a neighborhood $V$ of $J(G)$,
the identity map $i:(V,d_{s})\rightarrow
(V,d_{e})$ is a bi-Lipschitz map, where
$d_{s}$ and $d_{e}$ denote the
spherical and Euclidean distance,
respectively. In what follows, we often
use the above implicitly,
especially when we need to use the
facts in \cite{F} and \cite{Pe}.
\end{rem}
\section{Main Theorem A}
\label{Hausdimex}
In this section, we show Main Theorem A.
We investigate the estimate of the upper
box and Hausdorff dimensions
of Julia sets of expanding semigroups
using thermodynamic formalism in ergodic theory.
For the notation used in ergodic theory, see \cite{DGS} and \cite{W2}.
\begin{df}
\label{expandingdf}
Let $G=\langle f_{1},\cdots ,f_{m}\rangle$ be a finitely generated
rational semigroup. We say that $G$ is an {\bf expanding rational semigroup}
if $J(G)\neq \emptyset$ and
the skew product map $\tilde{f}: \Sigma _{m}\times \CCI \rightarrow
\Sigma _{m}\times \CCI$ associated with the generator system
$\{f_{1},\cdots ,f_{m}\}$ is expanding along fibers, i.e.,
there exists a positive constant $C$ and a constant $\lambda >1$
such that for each $n\in \NN$,
$$\inf\limits _{z\in \tilde{J}(\tilde{f})} \| (\tilde{f}^{n})'(z)\| \geq
C\lambda ^{n},$$
where we use $\| \cdot \|$ to denote the norm of the derivative with respect
to the spherical metric.
\end{df}
\begin{rem}
\label{exphyprem}
By Theorem 2.6, Theorem 2.8, and Remark 4 in \cite{S2},
we see that if
$G=\langle f_{1},\cdots ,f_{m}\rangle$ contains an
element of degree at least two, each
M\"{o}bius transformation in $G$ is neither the identity
nor an elliptic element, and $G$ is
{\bf hyperbolic}, i.e.,
the postcritical set $P(G)$ of $G$, which is defined
as:
$$ P(G):=\overline{\bigcup _{g\in G}\{\mbox{all critical values
of }g\}},$$
is included in $F(G)$,
then $G$ is expanding.
Conversely, if $G=\langle f_{1},\cdots ,f_{m}\rangle$
is expanding, then $G$ is hyperbolic and
each M\"{o}bius transformation in $G$ is loxodromic.
Hence,
the notion of expandingness does not depend on any choice of
a generator system for
a finitely generated rational semigroup.
\end{rem}
\begin{lem}
\label{elemexp}
Let $G=\langle f_{1},\cdots ,f_{m}\rangle$ be
a finitely generated expanding rational semigroup.
Suppose $\sharp J(G)\leq 2$. Then,
$\sharp J(G)=1$ and $J(G)$ is a common
repelling fixed point of any $f_{j}$.
\end{lem}
\begin{proof}
Suppose $\sharp J(G)=2$ and
let $J(G)=\{z_{1},z_{2}\}$.
Then, $f_{j}$ is a M\"{o}bius transformation,
for each $j=1,\cdots ,m$.
Since $G$ is expanding,
each $f_{j}$ is loxodromic.
We may assume that $z_{1}$ is a
repelling fixed point of $f_{1}$.
Then, since $f_{1}^{-1}(J(G))\subset J(G)$,
it follows that $z_{2}$ is an attracting fixed point
of $f_{1}$. This is a contradiction, however, since
$G$ is expanding.
Hence, $\sharp J(G)=1$.
\end{proof}

\begin{df}
Let $G$ be a rational semigroup and let 
$t$ be a non-negative number. We say
that a Borel probability measure $\tau$ on $\CCI$ is 
{\bf $t$-subconformal} (for $G$) if for
each $g\in G$ and for each Borel measurable set $A$ in $\CCI$,
$\tau (g(A))\leq \int _{A} \| g'\| ^{t}d\tau$.
Moreover, we set
\[s(G)=\inf \{t \mid \exists \tau : t \mbox{-subconformal measure}\} 
.\]
\end{df}
\begin{df}
Let $X$ be a compact metric space.
Let $f:X\rightarrow X$ be a continuous map:
\begin{enumerate}
\item We use $h(f)$ to denote the {\bf topological entropy}
of $f$ (see p83 in \cite{DGS}). We use $h_{\mu}(f)$ to denote the
{\bf metric entropy} of $f$ with respect to
an invariant Borel probability measure $\mu$ (see p60 in \cite{DGS}).
\item Furthermore, let $\varphi :X\rightarrow \RR$ be a
continuous function.
Then, we use $P(f,\varphi)$ to denote the
{\bf pressure} for the dynamics of $f$ and
the function $\varphi$ (see p141 in \cite{DGS}).
According to a well known fact: the variational principle
 (see p142 in \cite{DGS}),
we have
$$P(f,\varphi)=\sup \{h_{\mu}(f)+\int _{X}\varphi \ d\mu\},$$
where the supremum is taken over all $f$-invariant Borel probability
measures $\mu$ on $X$.
If an invariant probability measure $\mu$ attains the supremum
in this manner, then $\mu$ is called
an {\bf equilibrium state} for $(f, \varphi)$.
For more details on this notation and the variational principle,
see \cite{DGS} and \cite{W2}.
\item For a real-valued continuous function
$\varphi$ on $X$ and for each $n\in \NN$,
we define a continuous function
$S_{n}\varphi$ on $X$ as
$(S_{n}\varphi)(z)$ $=\sum _{j=0}^{n-1}\varphi (f^{j}(z))$.
Note that
$P(f^{n},S_{n}\varphi)=nP(f,\varphi)$(see Theorem 9.8 in \cite{W2}).
\end{enumerate}
\end{df}
\begin{df}Let $X$ be a compact metric space and
let $f:X\rightarrow X$ be a continuous map satisfying the fact that
there exists a number $k\in \NN$ such that
$\sharp f^{-1}(z)=k$ for each $z\in X$.
Let $\varphi$ be a continuous function on
$X$. We define an operator
$L=L_{\varphi}$ on $C(X)$ using
\[ L\psi (z)=\sum _{f(z')=z}
\exp (\varphi (z'))\psi (z').\] 
This is called
the {\bf transfer operator for $(f,\varphi)$}.
Note that $L_{\varphi}^{n}$ equals the transfer
operator for $(f^{n}, S_{n}\varphi)$,
for each $n\in \NN$.
\end{df}

\begin{lem}
\label{lem1}
Let $G= \langle f_{1},f_{2},\cdots f_{m} \rangle$
be a finitely generated expanding rational semigroup.
Let $\tilde{f}:\Sigma _{m}\times \CCI
\rightarrow \Sigma _{m}\times \CCI$ be the
skew product map associated with
$\{f_{1},\cdots ,f_{m}\}$.
Then, for each H\"{o}lder continuous
function $\varphi$ on $\tilde{J}(\tilde{f})$,
the transfer operator
$L_{\varphi}$ for $(\tilde{f}|_{\tilde{J}(\tilde{f})}, \varphi)$
on $C(\tilde{J}(\tilde{f}))$ satisfies the fact that
there exists a unique probability measure
$\tilde{\nu}=\tilde{\nu}_{\varphi}$ on $\tilde{J}(\tilde{f})$
satisfying all of the following:
\begin{enumerate}
\item \label{lem1-1} $L_{\varphi}^{\ast}\tilde{\nu}=\exp (P)
\tilde{\nu},$
where $P=P(\tilde{f}|_{\tilde{J}(\tilde{f})},
\varphi)$ is the pressure of $(\tilde{f}|_{\tilde{J}(\tilde{f})}, \varphi)$.
\item \label{lem1-2}
For each $\psi \in C(\tilde{J}(\tilde{f})),\ 
\| \frac{1}{(\exp (P))^{n}}L_{\varphi}^{n}\psi
-\tilde{\nu}(\psi)
\alpha _{\varphi}\| _{\tilde{J}(\tilde{f})} 
\rightarrow 0, n\rightarrow \infty$,
where we set $\alpha _{\varphi}=
\lim _{l\rightarrow \infty}
\frac{1}{(\exp (P))^{l}}
L_{\varphi}^{l}(1)\in C(\tilde{J}(\tilde{f}))$ and we use
$\| \cdot \| _{\tilde{J}(\tilde{f})}$ to denote the supremum norm on 
$\tilde{J}(\tilde{f})$.
\item \label{lem1-3}
$\alpha _{\varphi}\tilde{\nu}$ is $\tilde{f}$-invariant,
exact (hence ergodic) 
and is an equilibrium state for
$ (\tilde{f}| _{\tilde{J}(\tilde{f})}, \varphi)$.
\item \label{lem1-4}
$\alpha _{\varphi}(z)>0$ for each $z\in \tilde{J}(\tilde{f}).$
\end{enumerate}
\end{lem}
\begin{proof}
According to the Koebe distortion theorem and since the dynamics of
$\sigma :\Sigma _{m}\rightarrow \Sigma _{m}$ is expanding,
there exists a number $s\in \NN$ such that
the map $\tilde{f}^{s}:\tilde{J}(\tilde{f})\rightarrow
\tilde{J}(\tilde{f})$ satisfies condition I on page 123 in
\cite{W1} (each of $X_{0}, X$, and $\overline{X}$ in
\cite{W1} corresponds to $\tilde{J}(\tilde{f})$).
Furthermore,
by Proposition 3.2 (f) in \cite{S5}
and Lemma~\ref{elemexp},
we have the fact
that $\tilde{f}^{s}$ on $\tilde{J}(\tilde{f})$ satisfies
condition II on page 125 in \cite{W1}.
The map
$\mu \rightarrow L_{\varphi}^{\ast}
\mu /(L_{\varphi}^{\ast}\mu)(1)$ is
continuous on the space $M(\tilde{J}(\tilde{f}))$
of Borel probability measures on $\tilde{J}(\tilde{f}).$
Hence, this map has a fixed point $\tilde{\nu}$
based on the Schauder-Tychonoff fixed point theorem.
Let $\lambda =(L_{\varphi}^{\ast}\tilde{\nu})(1)$.
Then, $L_{\varphi}^{\ast}\tilde{\nu}
=\lambda \tilde{\nu}$. Hence, we have
$(L_{\varphi}^{s})^{\ast}\tilde{\nu}
=\lambda ^{s}\tilde{\nu}$.
By Theorem 8, Corollary 12, and the statement
on equilibrium states on page 140 in \cite{W1}, we get
$\lambda ^{s}=
\exp (P(\tilde{f}^{s}|_{\tilde{J}(\tilde{f})}, S_{s}\varphi))
=\exp(sP(\tilde{f}|_{\tilde{J}(\tilde{f})}, \varphi))$.
Hence, we obtain
$\lambda =\exp (P).$ 
The other
results also follow from
Theorem 8, Corollary 12, and the statement
on equilibrium states on page 140 in \cite{W1}.
\end{proof}

\noindent {\bf Notation:}
Let $G=\langle f_{1},\cdots ,f_{m}\rangle$ be a
finitely generated rational semigroup.
Let $\tilde{f}:\Sigma _{m}\times \CCI \rightarrow
\Sigma _{m}\times \CCI$ be the skew product
map associated with $\{f_{1},\cdots ,f_{m}\}$.
Suppose that no critical point of
$\tilde{f}$ exists in $\tilde{J}(\tilde{f})$. Then,
we define a function $\tilde{\varphi}$ on
$\tilde{J}(\tilde{f})$ as:
$\tilde{\varphi}((w,x)):=-\log \| (f_{w_{1}})'(x)\|$
for $(w,x)=((w_{1},w_{2},\cdots), x)\in \tilde{J}(\tilde{f})$.
\begin{lem}
\label{lem2}
Let $G= \langle f_{1},f_{2},\cdots f_{m} \rangle$ be a finitely generated
expanding rational semigroup.
Then, using the above notation, we have the following:
\begin{enumerate}
\item \label{lem2-1} The function
$P(t)=P(\tilde{f}|_{\tilde{J}(\tilde{f})}, t\tilde{\varphi})$
on $\RR$ is convex and strictly decreasing as
$t$ increases. Furthermore,
$P(t)\rightarrow -\infty$ as $t\rightarrow \infty$.
\item \label{lem2-2}
There exists a unique zero $\delta \geq 0$
of $P(t)$.
Furthermore, if $h(\tilde{f}|_{\tilde{J}(\tilde{f})})>0$
then $\delta >0$.
\item \label{lem2-3} There exists a unique probability
measure $\tilde{\nu}=\tilde{\nu}_{\delta \tilde{\varphi}}$
on $\tilde{J}(\tilde{f})$ such that
$M_{\delta}^{\ast}\tilde{\nu}=\tilde{\nu}$,
where $M_{\delta}$ is an operator on
$C(\tilde{J}(\tilde{f}))$ defined by
\begin{equation}
\label{mdelta}
M_{\delta}\psi ((w,x))=\sum _{\tilde{f} ((w',y))=(w,x)}
\frac{\psi ((w',y))}{\| (f_{w_{1}'})'(y)\| ^{\delta}}.
\end{equation}
Note that $M_{\delta}=L_{\delta \tilde{\varphi}}$.
\item \label{lem2-4}
$\delta$ satisfies the fact that
\begin{equation}
\label{deltaeq}
\delta =\frac{h_{\alpha \tilde{\nu}}(\tilde{f})}
{-\int _{\tilde{J}(\tilde{f})}\tilde{\varphi}\alpha d\tilde{\nu}}\leq
\frac{\log (\sum _{j=1}^{m}\deg (f_{j}))}
{-\int _{\tilde{J}(\tilde{f})}\tilde{\varphi} \alpha d\tilde{\nu}},
\end{equation}
where $\alpha =\lim _{l\rightarrow \infty}M_{\delta}^{l}(1)\in 
C(\tilde{J}(\tilde{f})).$
\end{enumerate}
\end{lem}
\begin{proof}
Using the variational principle, we have
$P(t)=\sup \{h_{\mu}(\tilde{f}|_{\tilde{J}(\tilde{f})}) +
\int _{\tilde{J}(\tilde{f})}t\tilde{\varphi}\ d\mu \}$,
where the supremum is taken over all $\tilde{f}$-invariant
Borel probability measures $\mu$ on $\tilde{J}(\tilde{f})$.
In addition, note that by Theorem 6.1 in \cite{S5}, we determine that
the topological entropy
$h(\tilde{f})$ of $\tilde{f}$
on $\Sigma _{m}\times \CCI$ is less than or equal
to $\log (\sum _{j=1}^{m}\deg (f_{j}))$.
By the variational principle:
$h(\tilde{f})=\sup \{h_{\mu}(\tilde{f})\mid
\tilde{f}_{\ast} \mu =\mu \}$
(see p138 in \cite{DGS} or Theorem 8.6 in \cite{W2}); it follows that
$$h_{\mu}(\tilde{f})\leq
\log (\sum _{j=1}^{m}(\deg (f_{j})))$$ 
for any
$\tilde{f}$-invariant
Borel probability measure $\mu$ on $\tilde{J}(\tilde{f})$.
Combining this with the fact that the dynamics of $\tilde{f}$ on 
$\tilde{J}(\tilde{f})$
is expanding, we see that the function $P(t)$ on $\RR$ is
convex, strictly decreasing as $t$ increases, and
$P(t)\rightarrow -\infty$ as $t\rightarrow \infty$.
Hence, there exists a unique number $\delta \in \RR$ satisfying
$P(\delta)=0$. Since $P(0)=h(\tilde{f}|_{\tilde{J}(\tilde{f})})$, we have $\delta >0$
if $h(\tilde{f}|_{\tilde{J}(\tilde{f})})>0$.
The statements \ref{lem2-3} and \ref{lem2-4} follow
from Lemma~\ref{lem1} and this argument.
\end{proof}
\begin{df}
We define an operator $\hat{M}_{\delta}$
acting on the space of all Borel measurable functions
on $\tilde{J}(\tilde{f})$ using the same formula as that for
$M_{\delta}$. (See (\ref{mdelta})).
\end{df}
We now show that $\hat{M}_{\delta}$ acts on
$L^{1}(\tilde{\nu})$ and that
$\hat{M}_{\delta}$ on $L^{1}(\tilde{\nu})$ is
a bounded operator, where $\tilde{\nu}=
\tilde{\nu}_{\delta \tilde{\varphi}}$.
\begin{lem}
\label{mdeltalem}
Let $G=\langle f_{1},\cdots ,f_{m}\rangle$ be a
finitely generated expanding rational semigroup.
Using the above notation,
we have the following:
\begin{enumerate}
\item \label{mdeltalem1}
Let $A$ be a Borel set in $\tilde{J}(\tilde{f})$.
If $\tilde{\nu}(A)=0$, then
$\tilde{\nu}(\tilde{f}^{-1}(A))=0$.
\item \label{mdeltalem2}
Let $\psi$ be a Borel measurable function
on $\tilde{J}(\tilde{f})$. Let $\{\psi _{n}\} _{n}$ be
a sequence of Borel measurable functions on
$\tilde{J}(\tilde{f})$. Suppose
$\psi _{n}(z)\rightarrow \psi (z)$ for almost
every $z\in \tilde{J}(\tilde{f})$ with respect to
$\tilde{\nu}$. Then,
we have $(\hat{M}_{\delta}\psi _{n})(z)\rightarrow
(\hat{M}_{\delta}\psi)(z)$ for almost every 
$z\in \tilde{J}(\tilde{f})$ with respect to
$\tilde{\nu}$.
\item \label{mdeltalem3}
If $\psi \in L^{1}(\tilde{\nu})$,
then $\hat{M}_{\delta}\psi \in L^{1}(\tilde{\nu})$.
Furthermore,
$\hat{M}_{\delta}$ is a bounded operator on
$L^{1}(\tilde{\nu})$ and the operator norm
$\| \hat{M}_{\delta}\|$ is equal to $1$.
\end{enumerate}
\end{lem}
\begin{proof}
Let $\mu =\alpha \tilde{\nu}$, where
$\alpha =\lim _{\l\rightarrow \infty}M_{\delta}^{l}1$.
Then, by
Lemma~\ref{lem1}-\ref{lem1-3}, we have
$\tilde{f}_{\ast}\mu =\mu$. Furthermore,
by Lemma~\ref{lem1}-\ref{lem1-4},
$\mu$ and $\tilde{\nu}$ are
absolutely continuous with respect to
each other. Hence, we obtain the statement
\ref{mdeltalem1}, and the statement \ref{mdeltalem2}
follows easily from this.

 We now show the statement \ref{mdeltalem3}.
First, we show the following claim:\\
Claim: for any $\psi \in C(\tilde{J}(\tilde{f}))$,
we have $\int |\hat{M}_{\delta}\psi |\ d\tilde{\nu}
\leq \int
|\psi |\ d\tilde{\nu}$.

 To show this claim, let $\psi \in C(\tilde{J}(\tilde{f}))$.
Let $\psi ^{+}=\max \{\psi ,0\}$ and
$\psi ^{-}=-\min \{\psi ,0\}$. Then,
we have $\psi =\psi ^{+}-\psi ^{-}$ and
$|\psi |=\psi ^{+}+\psi ^{-}$.
Since $M_{\delta}^{\ast}\tilde{\nu}=
\tilde{\nu}$, we obtain
$\int |\hat{M}_{\delta}\psi |\ d\tilde{\nu}
= \int|M_{\delta}\psi ^{+}-M_{\delta}\psi ^{-}|\
d\tilde{\nu}
\leq \int M_{\delta}\psi ^{+}\ d\tilde{\nu}
+\int M_{\delta}\psi ^{-}\ d\tilde{\nu}
=\int \psi ^{+}+\psi ^{-}\ d\tilde{\nu}
=\int |\psi |\ d\tilde{\nu}$. Hence, the above claim holds.

 Now,
let $\psi$ be a general element of
$L^{1}(\tilde{\nu})$.
Let $\{\psi _{n}\} _{n}$ be a sequence in
$C(\tilde{J}(\tilde{f}))$ such that
$\psi _{n}\rightarrow \psi$ in $L^{1}(\tilde{\nu})$.
We may assume that
$\psi _{n}(z)\rightarrow \psi (z)$ for almost
every $z\in \tilde{J}(\tilde{f})$ with respect to
$\tilde{\nu}$. Then, according to the statement \ref{mdeltalem2},
we have $(\hat{M}_{\delta}\psi _{n})(z)\rightarrow
(\hat{M}_{\delta}\psi)(z)$ for almost
every $z\in \tilde{J}(\tilde{f})$ with respect to
$\tilde{\nu}$. Using this claim,
$\{\hat{M}_{\delta}\psi _{n}\} _{n}$ is a
Cauchy sequence in $L^{1}(\tilde{\nu})$.
Hence, it follows that
$\hat{M}_{\delta}\psi \in L^{1}(\tilde{\nu})$.
Furthermore, we have
$\int |\hat{M}_{\delta}\psi |\ d\tilde{\nu}
=\lim_{n\rightarrow \infty}
\int |\hat{M}_{\delta}\psi _{n}|\ d\tilde{\nu}
\leq \lim _{n\rightarrow \infty}
\int |\psi _{n}|\ d\tilde{\nu}=
\int |\psi |\ d\tilde{\nu}$.
Hence, $\| \hat{M}_{\delta}\| \leq 1$. Since
$\int \hat{M}_{\delta}1\ d\tilde{\nu}=
\int 1\ d\tilde{\nu}=1$, we obtain
$\| \hat{M}_{\delta}\| =1$.
\end{proof}
We now show that the measure
$\tilde{\nu}=\tilde{\nu}_{\delta \tilde{\varphi}}$ is
``conformal''.
\begin{lem}
\label{nuconformal}
Let $G=\langle f_{1},\cdots ,f_{m}\rangle$ be a
finitely generated expanding rational semigroup.
Let $k\in \NN$
and let $A$ be a Borel set in $\tilde{J}(\tilde{f})$ such that
$\tilde{f}^{k}:A\rightarrow \tilde{f}^{k}(A)$ is
injective. Then, using the above notation,
we have
$\tilde{\nu}(\tilde{f}^{k}(A))=\int _{A}\| (\tilde{f}^{k})'\| 
^{\delta}\ d\tilde{\nu}$.
\end{lem}
\begin{proof}
We have $M_{\delta}^{k}\tilde{\nu}=\tilde{\nu}$ and
$M_{\delta}^{k}$ is a transfer operator
for $(\tilde{f}^{k}, \delta S_{k}\tilde{\varphi})$.
By Proposition 2.2 in \cite{DU} and
Lemma~\ref{mdeltalem}-\ref{mdeltalem3},
we obtain the statement.
\end{proof}
\begin{lem}
\label{nusubconf}
Let $G=\langle f_{1},\cdots ,f_{m}\rangle$ be a
finitely generated expanding rational semigroup.
Then, with our notation,
the probability measure
$\nu :=(\pi _{\CCI})_{\ast}(\tilde{\nu})$
is $\delta $-subconformal.
\end{lem}
\begin{proof}
First, note that
by Lemma~\ref{nuconformal}, it follows that for any
Borel set $B$ in $\Sigma _{m}\times \CCI$
we have $\tilde{\nu}(\tilde{f}^{k}(B))\leq
\sum _{j}\tilde{\nu}(\tilde{f}^{k}(B_{j}))=
\sum _{j}\int _{B_{j}}\| (\tilde{f}^{k})'\| ^{\delta}\
d\tilde{\nu}
= \int _{B}\| (\tilde{f}^{k})'\| ^{\delta}\ d\tilde{\nu}$,
where $B=\sum _{j}B_{j}$ is a measurable
partition such that $\tilde{f}^{k}|_{B_{j}}$ is
injective for each $j$.
Hence,
for any Borel set $A$ in $\CCI$
and any $w\in {\cal W}^{\ast}$ with $|w|=k$,
it follows that
$\nu (f_{w}(A))= \tilde{\nu}(\pi _{\CCI}^{-1}
(f_{w}(A)))=
\tilde{\nu}(\tilde{f}^{k}(\Sigma _{m}(w)\times A))
\leq \int _{\Sigma _{m}(w)\times A}\| (\tilde{f}^{k})'\| ^{\delta}
d\tilde{\nu}
\leq \int _{A}\| (f_{w})'\| ^{\delta}d\nu $.
\end{proof}

We now consider the Poincar\'{e} series and critical exponent for
a rational semigroup.
\begin{df}
Let $G$ be a rational semigroup.
We set $$A(G)=\overline{\cup _{g\in G}g(\{z\in \CCI
\mid \exists h\in G, h(z)=z, |h'(z)|<1\})}.$$
For any $s\geq 0$ and $x\in \CCI$,
we set $S(s,x)=\sum _{g\in G}\sum _{g(y)=x}$
$\|g'(y)\| ^{-s}$. Furthermore, we set
$S(x)=\inf \{s\geq 0\mid S(s,x)<\infty \}$
(If no $s$ exists with $S(s,x)<\infty$,
then we set $S(x)=\infty$). We set 
$s_{0}(G)=\inf \{ S(x)\mid x\in \CCI \} .$ 

 If $G$ is generated by finite elements
$\{f_{1},\cdots ,f_{m}\}$, then
for any $x\in \CCI$ and $t\geq 0$, we set
$T(t,x)=\sum _{w\in {\cal W}^{\ast}}
\sum _{f_{w}(y)=x}\| (f_{w})'(y)\| ^{-t}$
and $T(x)=\inf \{t\geq 0\mid
T(t,x)<\infty \}$
(If no $t$ exists with $T(t,x)<\infty$,
then we set $T(x)=\infty$).
Note that $S(x)\leq T(x)$.
\end{df}
\begin{lem}
\label{fsetn}
Let $G=\langle f_{1},\cdots
,f_{m}\rangle$ be a finitely
generated rational semigroup.
Let $\tilde{f}:\Sigma _{m}\times
\CCI \rightarrow \Sigma _{m}\times
\CCI$ be a skew product map
associated with $\{f_{1},\cdots ,f_{m}\}$.
Let $z\in \tilde{F}(\tilde{f})$ be a point.
Then, there exists a number $n\in \NN$ such that
$\pi _{\CCI}(\tilde{f}^{n}(z))\in F(G)$.
\end{lem}

\begin{proof}
Let $z\in \tilde{F}(\tilde{f})$ be a point.
Then, there exists a word $w\in {\cal W}^{\ast}$
and an open neighborhood $V$ of $\pi _{\CCI}(z)$
in $\CCI$ such that
$z\in \Sigma _{m}(w)\times V\subset \tilde{F}(\tilde{f}).$
Let $n=|w|$. Then,
$\tilde{F}(\tilde{f})\supset
\tilde{f}^{n}(\Sigma _{m}(w)\times V)
=\Sigma _{m}\times f_{w}(V)$.
Since $\pi _{\CCI}\tilde{J}(\tilde{f})=
J(G)$ (Proposition 3.2 in \cite{S5}), it follows
that $f_{w}(V)\subset F(G)$. Hence,
$\pi _{\CCI}\tilde{f}^{n}(z)=f_{w}(\pi _{\CCI}(z))
\in f_{w}(V)\subset F(G)$.
\end{proof}

\begin{lem}
\label{expconvj}
Let $G=\langle f_{1},\cdots
,f_{m}\rangle$ be a finitely
generated expanding rational semigroup.
Let $\tilde{f}:\Sigma _{m}\times
\CCI \rightarrow \Sigma _{m}\times
\CCI$ be a skew product map
associated with $\{f_{1},\cdots ,f_{m}\}$.
Let $z\in \Sigma _{m}\times \CCI$ be a
point with $x:=\pi _{\CCI}(z)\in \CCI
\setminus A(G)$. Then,
for each open neighborhood $V$ of $\tilde{J}(\tilde{f})$
in $\Sigma _{m}\times \CCI $, there exists a number
$l\in \NN$ such that $\cup _{n\geq l}(\tilde{f}^{n})^{-1}(z)
\subset V$. Furthermore, we have
$A(G)\cup P(G)\subset F(G)$ and if $x\in \CCI \setminus 
(A(G)\cup P(G))$,\ then 
$T(x)<\infty$.
\end{lem}

\begin{proof}
By Remark \ref{exphyprem}, we have $P(G)\subset F(G).$ 
Next we show $A(G)\subset F(G).$ 
Since $G$ is expanding, then using the Koebe
distortion theorem and 
$\pi _{\CCI }(\tilde{J}(\tilde{f}))=J(G)$ 
(Proposition 3.2 in \cite{S5}),\ we obtain
that there exist an $n\in \NN $ and a number 
$\varepsilon >0$ such that 
for each $a\in J(G)$ and 
each $w\in {\cal W}^{\ast }$ with 
$|w|=n$,\ we can take well-defined 
inverse branches of $f_{w}^{-1}$ on 
$B(a,\varepsilon )$ and any inverse branch 
$\gamma $ of $f_{w}^{-1}$ on $B(a,\varepsilon )$ 
satisfies $\gamma (B(a,\varepsilon ))\subset 
B(\gamma (a),\frac{1}{2}\varepsilon )$ and 
$\| \gamma '(y)\| \leq \frac{1}{2}$ for each 
$y\in B(a,\varepsilon ).$ Taking a small enough $\varepsilon $, 
it follows that  
for each $a\in J(G)$ and each $w\in {\cal W}^{\ast }$, 
we can take well-defined inverse branches $\gamma $ of 
$f_{w}^{-1}$ on $B(a,\varepsilon )$ and we have  
$$ \sup \{ \| \gamma '(y)\| \mid 
y\in B(a,\varepsilon ),\ a\in J(G),\ 
\gamma : \mbox{ a branch of } f_{w}^{-1},\ 
|w|=n\} \rightarrow 0$$ 
as $n\rightarrow \infty .$ 
Let $y\in \CCI $ be a point such that 
$g(y)=y$ and $|g'(y)|<1$ for some $g\in G.$ 
Suppose that there exist an element $h\in G$ 
and a point $a\in J(G)$ such that 
$h(y)\in B(a,\varepsilon ).$ Let 
$\gamma _{n}$ be a well-defined inverse branch 
of $(hg^{n})^{-1}$ on $B(a,\varepsilon )$ such that 
$\gamma _{n}(hg^{n}(y))=\gamma _{n}(h(y))=y.$ 
Then $|\gamma _{n}'(y)|\rightarrow \infty $ as 
$n\rightarrow \infty .$ This contradicts 
the previous argument. Hence 
$A(G)\subset \CCI \setminus B(J(G),\varepsilon )\subset F(G).$ 

 Next,\ suppose that there exists a sequence $(z_{j})$
in $\tilde{F}(\tilde{f})$ such that
$\tilde{f}^{n_{j}}(z_{j})=z$ and
$z_{j}\rightarrow z_{\infty}\in
\tilde{F}(\tilde{f})$ where
$n_{j}\in \NN$ with $n_{j}\rightarrow
\infty$ as $j\rightarrow \infty$.
Then, by Lemma~\ref{fsetn},
there exists a number $n\in \NN$ such that
$\pi _{\CCI}(\tilde{f}^{n}(z_{\infty}))\in F(G)$.
Let $x_{j}=\pi _{\CCI}(\tilde{f}^{n}(z_{j}))$
for each $j\in \NN$ and
let $x_{\infty}=\pi _{\CCI}\tilde{f}^{n}(z_{\infty})$.
Then, for each $j$ with $n_{j}>n$,
there exists an element $g_{j}\in G$ such that
$g_{j}(x_{j})=x$. Let
$\alpha =d(x,A(G))$. Since
$x_{j}\rightarrow x_{\infty}\in F(G)$,
we have
$\sharp \{j\mid d(g_{j}(x_{\infty}), x)
<\frac{\alpha }{2}\} =\infty$.

 By contrast, we have
$\sup \{d(f_{w}(x_{\infty}), A(G))\mid
|w|=n\} \rightarrow 0$ as $n\rightarrow \infty$.
For, if $P(G)\neq \emptyset$,
the above follows from Theorem 1.34 in \cite{S3}.
Even if $P(G)=\emptyset$, since $G$ is
expanding, by the Koebe distortion theorem,
then for each $z\in F(G)$,
$\overline{\cup _{g\in G}g(z)}\subset F(G)$. Using the
same argument as in the proof of Theorem 1.34 in
\cite{S3}, we obtain the above.

 Hence, we obtain a contradiction.
Therefore, we have shown that
for each open neighborhood $V$ of $\tilde{J}(\tilde{f})$
in $\Sigma _{m}\times \CCI $, there exists a number
$l\in \NN$ such that $\cup _{n\geq l}(\tilde{f}^{n})^{-1}(z)
\subset V$. 
If $x\in \CCI \setminus (A(G)\cup P(G))$,\ then 
since $G$ is expanding,
combining the above result with
the Koebe distortion theorem, we obtain
$T(x)<\infty$.
\end{proof}

\begin{df}
Let $E$ be a subset of $\CCI$,
$t\geq 0$ a number and $\beta >0$ a number.
We set
$$H_{\beta}^{t}(E):=
\inf \{\sum _{i=1}^{\infty}(\mbox{diam}(U_{i}))^{t}
\mid \mbox{diam}(U_{i})\leq \beta,
E\subset \cup _{i=1}^{\infty}U_{i}\}$$
and $H^{t}(E)=\lim _{\beta \rightarrow 0}H_{\beta}^{t}(E)$
with respect to the spherical metric on $\CCI$.
$H^{t}(E)$ is called the $t$-dimensional
(outer) Hausdorff
measure of $E$ with respect to the spherical metric.
Note that $H^{t}(E)$ is a Borel regular
measure on $\CCI $(see \cite{R}).
We set $\dim _{H}(E):=\sup \{t\geq 0\mid
H^{t}(E)=\infty \} =\inf \{t\geq 0\mid
H^{t}(E)=0\}$.
$\dim _{H}(E)$ is called the Hausdorff dimension of
$E$.
Furthermore, let $N_{r}(E)$ be the smallest
number of sets of spherical diameter $r$ that
can cover $E$.
We set $\overline{\dim}_{B}(E)=
\limsup\limits _{r\rightarrow 0}
\frac{\log N_{r}(E)}{-\log r}$.
$\overline{\dim}_{B}(E)$ is called
the upper box dimension of $E$.
\end{df}

\begin{lem}
\label{tsubconf}
Let $G=\langle f_{1},\cdots ,f_{m}\rangle$ be a
finitely generated expanding rational semigroup.
Let $\tau$ be a $t$-subconformal measure. Then,
there exists a positive constant $c$ such that
for each $r$ with $0<r<$ {\em diam} $\CCI$ and
each $x\in J(G)$, we have
$\tau (B(x,r))\geq cr^{t}$. Furthermore,
$H^{t}|_{J(G)}$ is absolutely continuous
with respect to $\tau ,$
$H^{t}(J(G))<\infty$ and
$\overline{\dim}_{B}(J(G))\leq t$.
\end{lem}
\begin{proof}
Let $\tau$ be a $t$-subconformal measure.
Using the argument in the proof of Theorem 3.4
in \cite{S2}, we find that
there exists a positive constant $c$ such that
for each $r$ with $0<r<$ diam $\CCI$ and
each $x\in J(G)$, $\tau (B(x,r))\geq cr^{t}$.
(Note that
for an estimate of this type,
we need only expandingness and we do not
need the strong open set condition used in \cite{S2}.)
By Proposition 2.2 in \cite{F},
we find that $H^{t}|_{J(G)}$ is
absolutely continuous with respect to
$\tau$. In particular,
$H^{t}(J(G))<\infty$. Furthermore,
by Theorem 7.1 in \cite{Pe}, we get
$\overline{\dim}_{B}(J(G))\leq t$.
\end{proof}
Using these arguments, we now demonstrate Main Theorem A.

\ 

\noindent {\bf Proof of Main Theorem A:} 
By Lemma~\ref{lem2}, we have 
that the function $P(t)$ has a unique 
zero $\delta $, there exists a unique probability 
measure $\tilde{\nu }$ on 
$\tilde{J}(f)$ such that $M_{\delta }^{\ast }\tilde{\nu }=\tilde{\nu },$  
and $\delta $ satisfies (\ref{deltaeq}).  

 By Lemma~\ref{tsubconf},
$\overline{\dim}_{B}(J(G))\leq s(G)$. 
By Lemma~\ref{expconvj}, we have 
$A(G)\cup P(G)\subset F(G).$ 
 Let $x\in \CCI \setminus (A(G)\cup P(G))$ be a point.
According to Theorem 4.2 in \cite{S2} and the fact that
$S(x)\leq T(x)<\infty$
(Lemma~\ref{expconvj}), we obtain
$s(G)\leq s_{0}(G)\leq S(x)\leq T(x)$.

 We show $\delta =T(x)$.
We consider the following two cases:\\
Case 1: $T(T(x),x)=\infty$.\\
Case 2: $T(T(x),x)<\infty$.\\
Suppose we have Case 1.
Let $z\in \Sigma _{m} \times \CCI$ be a
point with $\pi _{\CCI}(z)=x$.
Let $t_{n}$ be a sequence of real numbers
such that $t_{n}>T(x)$ for each $n\in \NN$ and
$t_{n}\rightarrow T(x)$.
For each $n\in \NN $, let $\mu _{n}$
be a Borel probability measure on
$\Sigma _{m}\times \CCI$ defined by:
$$\mu _{n}=\frac{1}{T(t_{n},x)}
\sum _{p\in \NN}\sum _{\tilde{f}^{p}(z')=z}
\| (\tilde{f}^{p})'(z')\| ^{-t_{n}}\delta _{z'}, $$
where $\delta _{z'}$ denotes the Dirac measure
concentrated at $z'$.
Since the space of
Borel probability measures on $\Sigma _{m}\times \CCI$
is compact, we may assume that
there exists a Borel probability measure
$\mu _{\infty}$ on $\Sigma _{m}\times \CCI$ such that
$ \mu_{n}\rightarrow \mu _{\infty}$ as $n\rightarrow \infty$,
with respect to the weak topology.
Then, by Lemma~\ref{expconvj},
we have supp $\mu _{\infty}\subset \tilde{J}(\tilde{f})$.
We now show the following claim:\\
Claim 1: For any Borel set $A$ in $\tilde{J}(\tilde{f})$
such that $\tilde{f}:A\rightarrow \tilde{f}(A)$ is
injective, we have
$\mu _{\infty}(\tilde{f}(A))=
\int _{A}\| (\tilde{f})'\| ^{T(x)}\ d\mu _{\infty}$.

 To show this claim,
let $A$ be a Borel set in $\Sigma _{m}\times \CCI$
such that $\tilde{f}:A\rightarrow \tilde{f}(A)$ is
injective. Then,
$\mu _{n}(\tilde{f}(A))=
\int _{A}\| (\tilde{f})'\| ^{t_{n}}\ d\mu _{n}
-\frac{1}{T(t_{n},x)}\sharp (\tilde{f}^{-1}(z)\cap A)$.
If $A$ satisfies that $\mu _{\infty}(\partial \tilde{f}(A))=
\mu _{\infty}(\partial A)=0$, then
letting $n\rightarrow \infty$ in the above,
it follows that
$\mu _{\infty}(\tilde{f}(A))=\int _{A}\| (\tilde{f})'\|
^{T(x)}\ d\mu _{\infty}$.

 Now let $B$ be a general Borel set in $\tilde{J}(\tilde{f})$
such that $\tilde{f}:B\rightarrow \tilde{f}(B)$ is injective.
Then, let $B=\sum _{j\in \NN}B_{j}$ be a
countable disjoint union of Borel sets $B_{j}$ satisfying the fact that
for each $j\in \NN$, there exists an open neighborhood
$W_{j}$ of $\overline{B_{j}}$ in $\Sigma _{m}\times \CCI$
such that $\tilde{f}:W_{j}\rightarrow \tilde{f}(W_{j})$ is
a homeomorphism.
Let $j$ be a fixed number and $K$ a fixed compact subset
of $W_{j}$. Then, for each $n\in \NN$, there exists a
number $\epsilon _{n}>0$ such that the
set $V_{n}:=\{z\in \Sigma _{m}\times \CCI \mid
d(z,K)<\epsilon _{n}\}$ satisfies
$V_{n}\subset W_{j},\ \mu _{\infty}(\partial V_{n})
=\mu _{\infty}(\partial (\tilde{f}(V_{n})))=0,\ 
\mu _{\infty}(\tilde{f}(V_{n})
\setminus \tilde{f}(K))<\frac{1}{n}$
and $\mu _{\infty}(V_{n}\setminus K)<\frac{1}{n}$.
For these sets $V_{n}$, by the previous argument, we have
$\mu _{\infty}(\tilde{f}(V_{n}))=\int _{V_{n}}
\| (\tilde{f})'\| ^{T(x)}\ d\mu _{\infty}$.
Letting $n\rightarrow \infty$, we obtain
$\mu _{\infty}(\tilde{f}(K))=
\int _{K}\| (\tilde{f})'\| ^{T(x)}\ d\mu _{\infty}$.
Next, for each $l\in \NN $, we can take a compact subset
$K_{l}$ of $B_{j}$ such that
$\mu _{\infty}(B_{j}\setminus K_{l})<\frac{1}{l}$ and
$\mu _{\infty}(\tilde{f}(B_{j})\setminus \tilde{f}(K_{l}))
<\frac{1}{l}$. For these sets $K_{l}$, using the above argument,
we have $\mu _{\infty}(\tilde{f}(K_{l}))=
\int _{K_{l}}\| (\tilde{f})'\| ^{T(x)}\ d\mu _{\infty}$.
Letting $l\rightarrow \infty $, we obtain
$\mu _{\infty}(\tilde{f}(B_{j}))=\int _{B_{j}}
\| (\tilde{f})'\| ^{T(x)}\ d\mu _{\infty}$.
Since $B=\sum _{j}B_{j}$ and $\tilde{f}$ is injective on
$B$, we obtain $\mu _{\infty}(\tilde{f}(B))=
\int _{B}\| (\tilde{f})'\| ^{T(x)}\ d\mu _{\infty}$.
Hence, we have shown Claim 1.

 Using Claim 1 and Proposition 2.2 in \cite{DU},
it follows that
$L_{T(x)\tilde{\varphi}}^{\ast}\mu _{\infty}
=\mu _{\infty}$.

 We now show that $\delta =T(x)$. Suppose $\delta <T(x)$.
Then, by Lemma~\ref{lem2}-\ref{lem2-1}, we have
$P(T(x))<0$. Then, for each
$\psi \in C(\tilde{J}(\tilde{f}))$,
we have
$ \mu _{\infty}(\psi)=(\exp P(T(x)))^{l}
\cdot \mu _{\infty}(\frac{L_{T(x)\tilde{\varphi}}\psi}
{(\exp (P(T(x))))^{l}})\rightarrow 0$ as $l\rightarrow \infty$,
by Lemma~\ref{lem1}-\ref{lem1-2}. Hence,
$\mu _{\infty}(\psi)=0$ and this implies a
contradiction. Suppose $T(x)<\delta$. Then, by a similar
argument to the one above, we get a contradiction.
Hence, $T(x)=\delta$.

 We now consider Case 2: $T(T(x),x)<\infty$.
Let $z\in \Sigma _{m}\times \CCI$ be a point with
$x=\pi _{\CCI}(z)$.
Then, we take Patterson's function(\cite{Pa}) $\Phi$:
i.e., $\Phi$ is a continuous, non-decreasing
function from $\RR _{+}:=\{t\in \RR \mid t\geq 0\}$ to
$\RR _{+}$ that satisfies the following:
\begin{enumerate}
\item $Q(t):=\sum _{n}\sum _{\tilde{f}^{n}(z')=z}
\Phi (\| (\tilde{f}^{n})'(z')\|)
\| (\tilde{f}^{n})'(z')\| ^{-t}$ converges for each
$t>T(x)$ and does not converge for each $t\leq T(x).$
\item For each $\epsilon >0$, there is a number
$r_{0}\in \RR _{+}$ such that
$\Phi (rs)\leq s^{\epsilon}\Phi (r)$ for each
$r>r_{0}$ and each $s>1$.
\end{enumerate}
Let $t_{n}$ be a sequence of $\RR$ such that
$t_{n}>T(x)$ for each $n\in \NN$,
$t_{n}\rightarrow T(x)$ as $n\rightarrow \infty$ and
the measures:
$$ \tau _{n}:=\frac{1}{Q(t_{n})}
\sum _{p}\sum _{\tilde{f}^{p}(z')=z}
\Phi (\| (\tilde{f}^{p})'(z')\|)
\| (\tilde{f}^{p})'(z')\| ^{-t_{n}}\delta _{z'}$$
tend to a Borel probability measure $\tau _{\infty}$
on $\Sigma _{m}\times \CCI$ as $n\rightarrow \infty$.
Then, by Lemma~\ref{expconvj},
we have supp $\tau _{\infty}\subset \tilde{J}(\tilde{f})$.
Furthermore, combining the argument in the proof of Claim 1
in Case 1 with the
properties of $\Phi$,
we find that for each Borel set $A$ in $\tilde{J}(\tilde{f})$
such that $\tilde{f}:A\rightarrow \tilde{f}(A)$ is injective,
$\tau _{\infty}(\tilde{f}(A))=
\int _{A}\| (\tilde{f})'\| ^{T(x)}\ d\tau _{\infty}$.
Combining this with the argument used in Case 1,
we obtain $\delta =T(x)$.

 Since $G$ is expanding and
$\nu$ is $\delta $-subconformal
(Lemma~\ref{nusubconf}), using an argument
in the proof of Theorem 4.4 in \cite{S2},
we obtain supp $\nu \supset J(G)$. Hence,
supp $\nu =J(G)$.

 Hence, we have shown Main Theorem A. 
\qed 

\begin{cor}
\label{cspmaincor}
Let $G= \langle f_{1},f_{2},\cdots f_{m} \rangle$
be a finitely generated
expanding rational semigroup.
Then,
$ \overline{\dim}_{B}(J(G))\leq \delta \leq \frac{\log (\sum _{j=1}^{m}\deg (f_{j}))}{\log 
\lambda}$,
where $\lambda$ denotes the number in Definition~\ref{expandingdf}.
\end{cor}

\begin{proof}
By Main Theorem A and (\ref{deltaeq}),
we have
\begin{align*}
\overline{dim}_{B}(J(G))
\leq \delta
& \leq
\frac{\log (\sum _{j=1}^{m}\deg (f_{j}))}
{-\int _{\tilde{J}(\tilde{f})}\tilde{\varphi}\alpha d\tilde{\nu}}\\
& = \frac{n\log (\sum _{j=1}^{m}\deg (f_{j}))}
{-\int _{\tilde{J}(\tilde{f})}S_{n}\tilde{\varphi}\ \alpha d\tilde{\nu}}\\
& \leq \frac{n\log (\sum _{j=1}^{m}\deg (f_{j}))}{\log C+n\log \lambda},
\end{align*}
for each $n\in \NN$. Letting $n\rightarrow \infty$,
we obtain the result.
\end{proof}

\section{Conformal measure}
In this section
we introduce the notion of ``conformal measure'',
which is needed in Main Theorem B.
\begin{df}
\label{cmdf}
\begin{enumerate}
\item \label{cmdf1}
Let $G$ be a rational semigroup.
Let $t\in \RR$ with $t \geq 0$.
We say that a Borel probability measure
$\tau$ on $J(G)$ is
$t$-{\bf conformal} (for $G$)
if
for any Borel set $A$ and
$g\in G,$ if $A,g(A)\subset J(G)$ and
$g:A\rightarrow g(A)$ is injective, then
$$\tau (g(A))=\int _{A}\| g'\| ^{t}\ d\tau .$$
\item \label{cmdf2}
Let $G=\langle f_{1},\cdots ,f_{m}\rangle$ be
a finitely generated rational semigroup.
We say that a Borel probability measure
$\mu$ on $J(G)$ satisfies the {\bf separating condition for
$\{f_{1},\cdots ,f_{m}\}$}
if
$\mu (f_{i}^{-1}(J(G))\cap f_{j}^{-1}(J(G)))=0$
for any $(i,j)$ with $i,j\in \{1,\cdots ,m\}$ and
$i\neq j$.
\end{enumerate}

\end{df}
We show some fundamental properties of
conformal measures.
\begin{lem}
\label{confsubconf}
Let $G=\langle f_{1},\cdots ,f_{m}\rangle$ be
a finitely generated rational semigroup.
Let $\tau$ be a $t$-conformal measure.
Then, $\tau$ is a $t$-subconformal measure.
\end{lem}
\begin{proof}
Let $A$ be a Borel set in $\CCI$ and
$g$ an element of $G$.
Let $J(G)=\sum B_{i}$ be a
measurable partition of $J(G)$
such that we can take the 
well-defined inverse branches of $g^{-1}$
on $B_{i}$, for each $i$
(we divide $J(G)$ into $\{B_{i}\}$ so that
for a critical value $c\in J(G)$ of $g$,
there exists an $i$ such that $B_{i}=\{c\}$).
Let $\{C_{i,j}\} _{j}$ be the
images of $B_{i}$ using the inverse branches of $g^{-1}$
so that $g:C_{i,j}\rightarrow B_{i}$ is
bijective for each $j$. Then, we have
$ \tau (g(A))=\tau (g(A)\cap J(G))=
\sum _{i}\tau (g(A)\cap B_{i})\leq
\sum _{i,j}\tau (g(A\cap C_{i,j}))
=\sum _{i,j}\int _{A\cap C_{i,j}}
\| g'\| ^{t}d\tau 
=\int _{A\cap \cup _{i,j}C_{i,j}}\| g'\| ^{t}d\tau 
\leq \int _{A}\| g'\| ^{t}d\tau .$
Hence, $\tau$ is $t$-subconformal.
\end{proof}

\begin{lem}
\label{cmlem1}
Let $G$ be
a rational semigroup.
Let $\tau$ be a Borel probability
measure on $J(G)$,
$g\in G$ an element,
and $V$ an open set in $\CCI$ with
$V\cap g^{-1}(J(G))\neq \emptyset$.
Suppose that $g:V\rightarrow g(V)$ is a
homeomorphism and
that for any Borel set $A$ in $V\cap g^{-1}(J(G))$,
$\tau (g(A))=\int _{A}\| g'\| ^{t}d\tau $.
Let $h:=(g|_{V})^{-1}:g(V)\rightarrow V$.
Then, we find that
for any Borel set $B$ in $g(V)\cap J(G)$,
$\tau (h(B))=\int _{B}\| h'\| ^{t} d\tau $.
\end{lem}
\begin{proof}
Let $\mu
:= h_{\ast}(\tau |_{g(V)\cap J(G)})$.
Then, by the assumption,
$d\mu =\| g'\| ^{t}d\tau '$,
where $\tau '=\tau |_{V\cap g^{-1}J(G)}$.
Let $B$ be a Borel set in
$g(V)\cap J(G)$. Then,
$\tau (h(B))=\int _{h(B)}\| g'\| ^{-t}\cdot
\| g'\| ^{t}d\tau $
$=\int _{h(B)}\| g'\| ^{-t}d\mu 
=\int _{B}\| g'\| ^{-t}\circ hd\tau 
=\int _{B}\| h'\| ^{t}d\tau $.
\end{proof}

\begin{lem}
\label{cmlemup}
Let $G=\langle f_{1},\cdots ,f_{m}\rangle$ be a
finitely generated rational semigroup.
Let $\tilde{f}:\Sigma _{m}\times \CCI
\rightarrow \Sigma _{m}\times \CCI$ be
the skew product map associated with
$\{f_{1},\cdots ,f_{m}\}$.
Let $\tilde{\tau}$ be a
Borel probability measure on $\tilde{J}(\tilde{f})$,
$n\in \NN$ an integer, and
$V$ an open set in $\Sigma _{m}\times \CCI$
such that $V\cap \tilde{J}(\tilde{f})\neq
\emptyset$.
Suppose that $\tilde{f}^{n}:
V\rightarrow \tilde{f}^{n}(V)$ is
a homeomorphism and that
for any Borel set $A$ in $\tilde{J}(\tilde{f})$,
$\tilde{\tau}(\tilde{f}^{n}(A))
=\int _{A}\|(\tilde{f}^{n})'\| ^{t}\ d\tilde{\tau}$.
Let $h=(\tilde{f}^{n}|_{V})^{-1}:\tilde{f}^{n}(V)
\rightarrow V$. Then,
we obtain the result that for any Borel set $B$
in $\tilde{f}^{n}(V)\cap \tilde{J}(\tilde{f})$,
$\tilde{\tau}(h(B))=
\int _{B}\| (\tilde{f}^{n})'(h)\| ^{-t}
\ d\tilde{\tau}$.
\end{lem}

\begin{proof}
This lemma can be shown using the same method as in the
proof of Lemma~\ref{cmlem1}.
\end{proof}
\begin{lem}
Let $G=\langle f_{1},\cdots ,f_{m}\rangle$ be
a finitely generated rational semigroup.
Let $\tau$ be a $t$-conformal measure satisfying
the separating condition for $\{f_{1},\cdots ,f_{m}\}$.
Suppose that for any $g\in G$,
if $c$ is a critical point of $g$ with
$g(c)\in J(G),$ then $\tau (\{c\})=0$.
Then, for any $k\in \NN$,
$\tau (f_{w}^{-1}(J(G))\cap f_{w'}^{-1}(J(G)))=0$
for any
$w=(w_{1},\cdots ,w_{k}), w'=(w'_{1},\cdots ,w'_{k})\in
\{1,\cdots ,m\} ^{k}$ with $w\neq w'$.
\end{lem}

\begin{proof}
Let $w=(w_{1},\cdots ,w_{k}), w'=(w'_{1},\cdots ,w'_{k})\in
\{1,\cdots ,m\} ^{k}$ with $w\neq w'$.
Let $1\leq u\leq k$ be the maximum
such that $w_{u}\neq w'_{u}$. If
$u=k$, then
$\tau (f_{w}^{-1}(J(G))\cap f_{w'}^{-1}(J(G)))
\leq \tau (f_{w_{k}}^{-1}(J(G))\cap
f_{w'_{k}}^{-1}(J(G)))=0$.
 Suppose that $u<k$.
Let $g=f_{w_{u+1}}\cdots f_{w_{k}}=
f_{w'_{u+1}}\cdots f_{w'_{k}}$.
Then, $f_{w}^{-1}(J(G))\cap f_{w'}^{-1}(J(G))
\subset g^{-1}(f_{w_{u}}^{-1}(J(G))\cap
f_{w'_{u}}^{-1}(J(G)))$. By
Lemma~\ref{cmlem1}, we have
$\tau (g^{-1}(f_{w_{u}}^{-1}(J(G))\cap
f_{w'_{u}}^{-1}(J(G))))=0$. Hence,
we obtain
$\tau (f_{w}^{-1}(J(G))\cap f_{w'}^{-1}(J(G)))=0$.
\end{proof}

\begin{df}
Let $G=\langle f_{1},\cdots ,f_{m}\rangle$ be a
rational semigroup. Suppose that
for each $g\in G$, no
critical value of $g$ exists in $J(G)$. Let
$t\in \RR$.
We define an operator
$N_{t}:C(J(G))\rightarrow C(J(G))$ as follows: \\
$(N_{t}\psi)(z)=\sum_{j=1}^{m}\sum_{f_{j}(y)=z}
\| f_{j}'(y)\| ^{-t}\psi (y)$ for each
$\psi \in C(J(G))$.
\end{df}

\begin{lem}
\label{ntlem}
Let $G=\langle f_{1},\cdots ,f_{m}\rangle$ be a
rational semigroup. Suppose that
for each $g\in G$, no
critical value of $g$ exists in $J(G)$.
Let $\tilde{f}:\Sigma _{m}\times \CCI \rightarrow
\Sigma _{m}\times \CCI$ be the skew
product map associated with
$\{f_{1},\cdots ,f_{m}\}$.
Then, we have the following commutative diagram:
$$
\begin{CD}
C(J(G)) @>{N_{t}}>>C(J(G))\\
@V{(\pi _{\CCI})^{\ast}}VV
@VV{(\pi _{\CCI})^{\ast}}V\\
C(\tilde{J}(\tilde{f})) @ >>{L_{t\tilde{\varphi}}}>
C(\tilde{J}(\tilde{f})).
\end{CD}
$$
\end{lem}

\begin{proof}
Let $\psi \in C(\tilde{J}(\tilde{f}))$ and
$(w,x)\in \tilde{J}(\tilde{f})$. Then,
$((\pi _{\CCI})^{\ast}N_{t}\psi)((w,x))
=(N_{t}\psi)(x)$ $=
\sum_{j=1}^{m}\sum_{f_{j}(y)=x}
\| f_{j}'(y)\| ^{-t}\psi (y)$.
Conversely,
$(L_{t\tilde{\varphi}}(\pi _{\CCI})^{\ast}\psi)((w,x))
=\sum _{\tilde{f}((w',y))=(w,x)}
\| f_{w'_{1}}'(y)\| ^{-t}((\pi _{\CCI})^{\ast}
\psi)((w',y))=
\sum _{j=1}^{m}\sum _{f_{j}(y)=x}\| f_{j}'(y)\| ^{-t}
\psi (y)$.
\end{proof}

\begin{lem}
\label{cmnt}
Let $G=\langle f_{1},\cdots ,f_{m}\rangle$ be a
rational semigroup. Suppose that
for each $g\in G$, no
critical value of $g$ exists in $J(G)$. Then,
we have the following:
\begin{enumerate}
\item \label{cmnt1}
Let $\tau$ be a $t$-conformal measure.
Then, we have $N_{t}^{\ast}\tau \geq \tau$; i.e.,
for each $\psi \in C(J(G))$ such that
$0\leq \psi (z)$ for each $z\in J(G)$,
we have $(N_{t}^{\ast}\tau)(\psi)\geq \tau (\psi)$.
\item \label{cmnt2}
If $\tau$ is a $t$-conformal measure satisfying
the separating condition for \\
$\{f_{1},\cdots ,f_{m}\}$,
then $N_{t}^{\ast}\tau =\tau$.
\item \label{cmnt3}
If $\tau$ is a $t$-conformal measure
satisfying $N_{t}^{\ast}\tau =\tau $, then
$\tau$ satisfies the separating condition for
$\{f_{1},\cdots ,f_{m}\}$.
\end{enumerate}
\end{lem}

\begin{proof}
Let $J(G)=\sum _{i=1}^{u}B_{i}$ be a
measurable partition of $J(G)$ such that
for each $j=1,\cdots ,m$ and $i=1,\cdots ,u$,
we can take the well-defined inverse branches
of $f_{j}^{-1}$ on $B_{i}$.
Then, for any Borel probability measure
$\tau$ on $J(G)$ and any
$\psi \in C(J(G))$,
we have
\begin{align*}
\int N_{t}\psi \ d \tau
& =
\int \sum_{j=1}^{m}\sum _{f_{j}(y)=z}\|
f_{j}'(y)\| ^{-t}\psi (y)\ d\tau (z)\\
& =\sum _{j}\sum _{i}\sum _{\gamma}
\int _{B_{i}}\|f_{j}'(\gamma (z))\| ^{-t}
\psi (\gamma (z))\ d\tau (z),
\end{align*}
where $\gamma$ runs over all inverse branches
of $f_{j}^{-1}$ on $B_{i}$.
Suppose that $\tau$ is $t$-conformal.
Then, we have
\begin{align*}
\int _{B_{i}}\| f_{j}'(\gamma (z))\| ^{-t}
\psi (\gamma (z))\ d\tau (z)
& = \int _{\gamma (B_{i})}
\| f_{j}'(x)\| ^{-t}\psi (x)\
d(\gamma _{\ast}(\tau |_{B_{i}}))(x)\\
& = \int _{\gamma (B_{i})}
\| f_{j}'(x)\| ^{-t}\psi (x)\cdot
\| f_{j}'(x)\| ^{t}\ d\tau (x)\\
& = \int _{\gamma (B_{i})}\psi (x)\ d\tau (x).
\end{align*}
Hence,
$\int N_{t}\psi \ d\tau =
\sum _{j}\sum _{i}\sum _{\gamma}
\int _{\gamma (B_{i})}\psi \ d\tau  $,
which is larger than or equal to
$\int _{J(G)}\psi \ d\tau $ if
$0\leq \psi (z)$ for each $z\in J(G)$,
since $J(G)=\cup _{j=1}^{m}f_{j}^{-1}(J(G))$
(Lemma 1.1.4 in \cite{S1}).
Furthermore, if $\tau$ is a $t$-conformal measure
satisfying the separating condition for
$\{f_{1},\cdots ,f_{m}\}$, then
for each $\psi \in C(J(G))$, we have
$\int N_{t}\psi \ d\tau
=\sum _{j}\sum _{i}\sum _{\gamma}\int _{\gamma (B_{i})}
\psi \ d\tau =\int _{J(G)}\psi \ d\tau $,
by $J(G)=\cup _{j=1}^{m}f_{j}^{-1}(J(G))$.

 We now show the statement \ref{cmnt3}.
Let $\tau$ be a $t$-conformal measure
satisfying $N_{t}^{\ast}\tau =\tau$.
Let $\psi \in C(J(G))$ be an element with
$\psi (x)\geq 0$ for each $x\in J(G)$. Then,
by the above argument, it follows that 
$\int _{J(G)}N_{t}\psi \ d\tau =
\sum _{j}\sum _{i}\sum _{\gamma}\int _{\gamma (B_{i})}
\psi \ d\tau \geq
\int _{J(G)}\psi \ d\tau $, where
$\gamma$ runs over all inverse branches of $f_{j}^{-1}$
on $B_{i}$. Since
$N_{t}^{\ast}\tau =\tau$, we have the
equality shown above. Hence, $\tau$ satisfies
the separating condition for $\{f_{1},\cdots ,f_{m}\}$.
\end{proof}

\begin{lem}
\label{convnt}
Let $G=\langle f_{1},\cdots ,f_{m}\rangle$ be
a finitely generated expanding rational semigroup.
Let $\delta$ be the number in
Lemma~\ref{lem2}, $t\geq 0$ a number, and
$\tilde{\nu}_{t\tilde{\varphi}}$ 
the Borel probability measure
on $\tilde{J}(\tilde{f})$ that is obtained in
Lemma~\ref{lem1}(
the unique fixed point of
$L_{t\tilde{\varphi}}^{\ast}$).
Let $\nu _{t}:=(\pi _{\CCI})_{\ast}
\tilde{\nu}_{t\tilde{\varphi}}$.
Then, we have the following:
\begin{enumerate}
\item \label{convnt1}
$\nu :=\nu _{\delta}$ satisfies
$N_{\delta}^{\ast}\nu =\nu$.
\item \label{convnt2}
$\frac{1}{(\exp (P(t)))^{l}}N_{t}^{l}\psi
\rightarrow
\nu _{t}(\psi)\cdot
\lim\limits _{l\rightarrow \infty}
N_{t}^{l}1$ in
$C(J(G))$, where
$P(t)=P(\tilde{f}|_{\tilde{J}(\tilde{f})}, t\tilde{\varphi})$.
\item \label{convnt3}
If $\tau$ is a
Borel probability measure
on $J(G)$ such that
$N_{t}^{\ast}\tau =\tau$,
then $t=\delta$ and $\tau =\nu$.
\end{enumerate}
\end{lem}

\begin{proof}
By Lemma~\ref{ntlem},
we obtain the statement \ref{convnt1}.
Since $\pi _{\CCI}(\tilde{J}(\tilde{f}))=J(G)$
(Proposition 3.2
in \cite{S5}), we find that
$(\pi _{\CCI})^{\ast}:C(J(G))\rightarrow
C(\tilde{J}(\tilde{f}))$ is an isometry
with respect to the supremum norms.
Hence, by Lemma~\ref{lem1} and Lemma~\ref{ntlem},
we find that
$\{\frac{1}{(\exp (P(t)))^{l}}N_{t}^{l}\psi\}
_{l\in \NN}$ is a Cauchy sequence in
$C(J(G))$. Let
$\psi _{0}=\lim _{l\rightarrow \infty}
\frac{1}{(\exp (P(t)))^{l}}N_{t}^{l}\psi$.
Then, by Lemma~\ref{lem1}, we obtain
\begin{align*}
(\pi _{\CCI})^{\ast}\psi _{0}
& = \lim _{l\rightarrow \infty}\frac{1}{(\exp (P(t)))^{l}}
L_{t\tilde{\varphi}}^{l}
(\pi _{\CCI})^{\ast}\psi \\
& = \tilde{\nu}_{t\tilde{\varphi}}
((\pi _{\CCI})^{\ast}\psi)\cdot \alpha _{t\tilde{\varphi}}\\
& = \nu _{t}(\psi)\cdot \lim _{\l\rightarrow \infty}
L_{t\tilde{\varphi}}^{l}(\pi _{\CCI})^{\ast}1\\
& = (\pi _{\CCI})^{\ast}(\nu _{t}(\psi)\cdot
\lim _{l\rightarrow \infty}N_{t}^{l}1).
\end{align*}
Hence, we obtain
$\psi _{0}=\nu _{t}(\psi)\cdot
\lim _{l\rightarrow \infty}N_{t}^{l}1$.

 Now, let $\tau$ be a Borel probability measure
on $J(G)$ such that $N_{t}^{\ast}\tau =\tau$.
Then for any $\psi \in C(J(G))$, we have
$\tau (\psi)=((N_{t}^{l})^{\ast}\tau)(\psi)
=\tau (N_{t}^{l}\psi)=
((\exp (P(t)))^{l})\cdot \tau (\frac{1}{(\exp (P(t)))^{l}}
N_{t}^{l}\psi)$ for any $l\in \NN$.
Since $\frac{1}{(\exp (P(t)))^{l}}
N_{t}^{l}\psi \rightarrow \nu _{t}(\psi)\cdot
\lim _{l\rightarrow \infty}N_{t}^{l}1$ and
$N_{t}^{\ast}\tau =\tau$,
we have $\tau (\frac{1}{(\exp (P(t)))^{l}}
N_{t}^{l}\psi)\rightarrow \nu _{t}(\psi)$ as
$l\rightarrow \infty$. Hence, it must be true that
$P(t)=0$, otherwise, 
we have
$\tau (\psi)=0$ for all $\psi$ or
$\tau (\psi)$ is not bounded,
both of which produce a contradiction.
Hence, it follows that $t=\delta$.
Further, by the above argument,
we obtain $\tau (\psi)=
\nu _{\delta}(\psi)$ for any $\psi \in C(J(G))$.
\end{proof}

\begin{lem}
\label{uniquecm}
Let
$G= \langle f_{1},f_{2},\cdots f_{m} \rangle$
be a finitely generated expanding rational semigroup.
Then, under the notation in
Lemma~\ref{convnt}, we have the following:
\begin{enumerate}
\item \label{uniquecm1}
If there exists a $t$-conformal measure
$\tau $, then $s(G)\leq t\leq \delta$.
\item \label{uniquecm2}
If there exists a $t$-conformal measure $\tau$ satisfying
the separating condition for $\{f_{1},\cdots ,f_{m}\}$,
then $t=\delta$ and $\tau =\nu$.
\end{enumerate}
\end{lem}

\begin{proof}
First, we show the statement \ref{uniquecm1}.
By Lemma~\ref{cmnt},
we have $N_{t}^{\ast}\tau \geq \tau$.
Hence, for each $\psi \in C(J(G))$ such that
$0\leq \psi (z)$ for each $z\in J(G)$, we have
$\tau (N_{t}^{l}\psi)\geq \tau (\psi)$ for each
$l\in \NN$. Suppose that $t>\delta$. Then,
by Lemma~\ref{lem2}-\ref{lem2-1},
$P(t)<0$. Hence, we obtain
$\tau (N_{t}^{l}\psi)=(\exp (P(t)))^{l}\cdot
\tau (\frac{N_{t}^{l}\tau}{(\exp(P(t))^{l}})\rightarrow
0$ as $l\rightarrow \infty $, by Lemma~\ref{convnt}-\ref{convnt2}.
Hence, $\tau (\psi)=0$ for each
$\psi \in C(J(G))$ such that
$0\leq \psi (z)$ for each $z\in J(G)$. This is a
contradiction, since $\tau (1)=1$. Hence,
$t\leq \delta$ must hold. By Lemma~\ref{confsubconf}, we
have $s(G)\leq t$. Hence, the statement \ref{uniquecm1} holds.

 Next, we show the statement \ref{uniquecm2}.
By Lemma~\ref{cmnt}-\ref{cmnt2},
we have $N_{t}^{\ast}\tau =\tau$.
Hence, by Lemma~\ref{convnt}-\ref{convnt3}, it follows that
$t=\delta$ and $\tau =\nu$.
\end{proof}

\begin{lem}
\label{hausconf}
Let $G$ be a
rational semigroup
and $t\geq 0$ a number.
Suppose that
$0<H^{t}(J(G))<\infty$.
Let $\tau =\frac{H^{t}|_{J(G)}}{H^{t}(J(G))}$.
Then, $\tau$ is a $t$-conformal measure.
\end{lem}
\begin{proof}
Suppose that $t=0$. Then,
each point $z\in \CCI$ satisfies
$H^{0}(\{z\})=1$. Since
we assume $0<H^{t}(J(G))<\infty$,
it follows that $1\leq \sharp (J(G))<\infty$.
Then, $G$ consists of degree 1 maps and
it is easy to see that $\tau$ is
$0$-conformal.

 Suppose that $t>0$.
Then, $H^{t}$ has no point mass.
Let $g\in G$ be an element.
 
 Step 1: For a critical point $c$ of $g$ in $J(G)$,
we have
$0=\tau (g(\{c\}))=\int _{\{c\}}\| g'\| ^{t}\
d\tau $.

 Step 2: Let $W$ be a non-empty open
set in $\CCI$ such that
$g:W\rightarrow g(W)$ is a diffeomorphism.
Let $K$ be a compact subset of
$W$ and $c>0$ a number.
Let $A$ be a Borel set
such that
$A\subset \{z\in g^{-1}(J(G))\mid d(z,A)<c\}
\subset K\cap g^{-1}(J(G))$.
Then, we show
the following claim:\\
Claim 1: we have
$H^{t}(g(A))=\int _{A}\| g'\| ^{t}\ dH^{t}$.

 To show this claim, let $\epsilon >0$ be a
given number. Let
$K=\sum _{i=1}^{l}K_{i}$ be a
disjoint union of Borel sets
$K_{i},\ \{z_{i}\} _{i=1}^{l}$ a
set with $z_{i}\in K_{i}$ for each $i$,
and $\xi >0$ a real number, such that:
\begin{enumerate}
\item $1-\epsilon \leq
\frac{\| g'(z)\|}{\| g'(z_{i})\|}
\leq 1+\epsilon $, for each $z\in B(K_{i}, \xi)$, and
\item
$(1-\epsilon)\| g'(z_{i})\|$ diam
$C\leq$ diam $g(C)\leq (1+\epsilon)\| g'(z_{i})\|$
diam $C$, for each subset $C$ of
$B(K_{i}, \xi)$.
\end{enumerate}
Let $i\in \NN \ (1\leq i\leq l)$ be a fixed number.
Let $\beta$ be a number with $0<\beta <\xi$.
Let $\{U_{p}\} _{p=1}^{\infty}$ be a
sequence of sets such that
$A\cap K_{i}\subset
\cup _{p=1}^{\infty}U_{p},
A\cap K_{i}\cap U_{p}\neq \emptyset$ for each $p\in \NN$ and
diam $U_{p}\leq \beta$ for each $p\in \NN$.
Then, since $\beta <\xi$, we have
$U_{p}\subset B(K_{i}, \xi)$ for each $p\in \NN$.
Hence, we obtain
$g(A\cap K_{i})\subset \cup _{p=1}^{\infty}
g(U_{p})$,
diam $g(U_{p})\leq (1+\epsilon)
\| g'(z_{i})\|$ diam $U_{p}$ for each $p\in \NN$
and
$\sum _{p=1}^{\infty}($
diam $g(U_{p}))^{t}\leq
(1+\epsilon)^{t}\| g'(z_{i})\| ^{t}
\sum _{p=1}^{\infty}($ diam $U_{p})^{t}$.
This implies that
$H^{t}_{(1+\epsilon)\| g'(z_{i})\| \beta}
(g(A\cap K_{i}))\leq (1+\epsilon)^{t}
\| g'(z_{i})\| ^{t}\sum _{p=1}^{\infty}
($ diam $U_{p})^{t}$.
Hence, we obtain
$H^{t}_{(1+\epsilon)\| g'(z_{i})\| \beta}
(g(A\cap K_{i}))\leq (1+\epsilon)^{t}
\| g'(z_{i})\| ^{t}H_{\beta}^{t}(A\cap K_{i})$.
Then, we obtain
$H^{t}(g(A\cap K_{i}))$ $\leq (1+\epsilon)^{t}
\| g'(z_{i})\| ^{t}\cdot H^{t}(A\cap K_{i})$,
letting $\beta \rightarrow 0$.
Similarly, we obtain
$H^{t}(A\cap K_{i})\leq
(1-\epsilon)^{-t}\| g'(z_{i})\| ^{-t}\cdot$
$H^{t}(g(A\cap K_{i}))$.
Hence, it follows that 
$(1-\epsilon)^{t}\| g'(z_{i})\| ^{t}
H^{t}(A\cap K_{i})\leq
H^{t}(g(A\cap K_{i}))\leq
(1+\epsilon)^{t}\| g'(z_{i})\| ^{t}
H^{t}(A\cap K_{i})$.
Moreover,
$(1-\epsilon)^{t}\| g'(z_{i})\| ^{t}\cdot$
$H^{t}(A\cap K_{i})\leq
\int _{A\cap K_{i}}\| g'\| ^{t}\
dH^{t}\leq
(1+\epsilon)^{t}\| g'(z_{i})\| ^{t}
H^{t}(A\cap K_{i})$.
Hence, we obtain
$$| H^{t}(g(A\cap K_{i}))-\int _{A\cap K_{i}}\| g'\| ^{t}\ dH^{t}|
\leq ((1+\epsilon)^{t}-(1-\epsilon)^{t})
\| g'(z_{i})\| ^{t} H^{t}(A\cap K_{i}).$$
This implies that
$|H^{t}(g(A))-\int _{A}\| g'\| ^{t}\
dH^{t}|
\leq ((1+\epsilon)^{t}-(1-\epsilon)^{t})\cdot
\max\limits _{z\in K}\| g'(z)\| ^{t}\cdot
\sum _{i=1}^{l}H^{t}(A\cap K_{i})$.
Since this inequality holds for each
$\epsilon >0$,
it follows that $H^{t}(g(A))=
\int _{A}\| g'\| ^{t}\ dH^{t}$.
Hence, we have shown Claim 1.

 Step 3: Let $B$ be a general Borel subset
of $g^{-1}(J(G))$ such that
$g:B\rightarrow g(B)$ is injective.
Let $B=\sum _{u=1}^{q}\{c_{u} \} \amalg
\sum _{v=1}^{\infty}B_{v}$
be a disjoint union of
Borel sets such that
each $c_{u}$ is a critical point of
$g$ (if one exists) and
for each $B_{v}$
there exists an open set $W_{v}$ in
$\CCI$ such that
$\overline{B_{v}}\subset W_{v}$ and
$g:W_{v}\rightarrow g(W_{v})$ is a
diffeomorphism.
Then, by Steps 1 and 2,
we obtain
$0=\tau (\{g(c_{u})\})
=\int _{\{c_{u}\}}\| g'\| ^{t}
\ d\tau $ for each $u$, and
$\tau (g(B_{v}))=
\int _{B_{v}}\| g'\| ^{t}\ d\tau $ for
each $v$. Combining this result with
the fact that $g:B\rightarrow g(B)$ is
injective,
it follows that
$\tau (g(B))=\sum _{u=1}^{q}\tau (\{g(c_{u})\})
+$ $\sum _{v=1}^{\infty}\tau (g(B_{v}))=
\sum _{v=1}^{\infty}
\int _{B_{v}}\| g'\| ^{t}\ d\tau $
$=\int _{B}\| g'\| ^{t}\ d\tau $.

 Hence, we have shown Lemma~\ref{hausconf}.
\end{proof}

\begin{lem}
\label{nowhere}
Let $G$ be a rational semigroup.
Let $\tau$ be a $t$-subconformal measure
for some $t\in \RR$. Suppose that
supp $\tau =J(G)$. Let $g\in G$.
Then, each Borel subset $A$ of
$g^{-1}(J(G))$ with $\tau (A)=0$ has
no interior points with respect to
the induced topology on $g^{-1}(J(G))$.

\end{lem}

\begin{proof}
Suppose there exists an open set $U$ of $\CCI$ such that
$A\supset U\cap g^{-1}(J(G))\neq \emptyset$.
Then, it follows that
$\tau (g(U))=\tau (g(U)\cap J(G))
=\tau (g(U\cap g^{-1}(J(G))))$
$\leq \int _{U\cap g^{-1}(J(G))}\| g'\| ^{t}\ d\tau 
=0$. This is a contradiction because
we assume supp $\tau =J(G)$.
\end{proof}
The following proposition is needed to show
Main Theorem B.
\begin{prop}
\label{deltasep}
Let $G=\langle f_{1},\cdots ,f_{m}\rangle$ be a
finitely generated expanding rational semigroup.
Let $\tilde{f}:\Sigma _{m}\times \CCI
\rightarrow \Sigma _{m}\times \CCI$ be the
skew product map associated with
$\{f_{1},\cdots ,f_{m}\}$.
Let $\delta$ be a number in Lemma~\ref{lem2}.
Let $\nu :=(\pi _{\CCI})_{\ast}
(\tilde{\nu}_{\delta \tilde{\varphi}})$.
Suppose that $0<H^{\delta}(J(G))$. Then,
we have $H^{\delta}(J(G))<\infty$,
$\nu =\frac{H^{\delta}|_{J(G)}}{H^{\delta}(J(G))}$, and
$\nu$ is a $\delta$-conformal measure
satisfying the separating condition with
respect to $\{f_{1},\cdots ,f_{m}\}$.
Furthermore,
$f_{i}^{-1}(J(G))\cap f_{j}^{-1}(J(G))$
is nowhere dense in
$f_{j}^{-1}(J(G))$, for each
$(i,j)$ with $i\neq j$.
\end{prop}
\begin{proof}
By Lemma~\ref{nusubconf} and Lemma~\ref{tsubconf},
we obtain that $\nu$ is a $\delta $-subconformal measure,
$H^{\delta}|_{J(G)}$ is absolutely continuous
with respect to $\nu$, and $H^{t}(J(G))<\infty$.
Let $\tau :=\frac{H^{\delta}|_{J(G)}}{H^{\delta}(J(G))}$.
Let $\varphi \in L^{1}(\nu)$ be the density function
such that $\tau (A)=\int _{A}\varphi \ d\nu$ for
any Borel subset $A$ of $J(G)$.
We show the following claim:\\
Claim 1: We have $(\varphi \circ \pi _{\CCI}
\circ \tilde{f})(z)\geq (\varphi \circ \pi _{\CCI})(z)$
for almost every $z\in \tilde{J}(\tilde{f})$ with respect
to $\tilde{\nu}:=\tilde{\nu}_{\delta \varphi}$.

 To show this claim,
let $j\ (1\leq j\leq m)$ be a number and
$A$ an open subset of $J(G)$ such that
we can take a well-defined inverse branch
$\gamma$ of $f_{j}^{-1}$ on $A$.
By Lemma~\ref{hausconf},
$\tau$ is $\delta $-conformal.
Hence, for each Borel subset
$B$ of $A$, we have
$\tau (B)=\int _{\gamma (B)}\| f_{j}'\| ^{\delta}
\ d\tau 
=\int _{\gamma (B)}\| f_{j}'\| ^{\delta}
\varphi \ d\nu $.
Moreover,
by Lemma~\ref{nusubconf}, we have
$\nu$ is $\delta $-subconformal.
Hence, we obtain
$\tau (B)=
\int _{B}\varphi \ d\nu $
$= \int _{A}(\varphi \circ f_{j}\circ \gamma)
\cdot (1_{\gamma (B)}\circ \gamma)\ d\nu
=\int _{A}(\varphi \circ f_{j})\cdot
1_{\gamma (B)}\
d(\gamma _{\ast}(\nu |_{A}))
\leq \int _{\gamma (B)}
\| f_{j}'\| ^{\delta}((\varphi \circ
f_{j}))\ d\nu $.
Hence, we obtain
$\varphi (x)\leq (\varphi \circ f_{j})(x)$
for almost every $x\in \gamma (A)$ with respect to
$\nu$.
It follows that
for each $j=1,\cdots ,m$, we have
$\varphi (x)\leq (\varphi \circ f_{j})(x)$ for
almost every $x\in f_{j}^{-1}(J(G))$ with respect to
$\nu$. This implies that
for each $j=1,\cdots ,m$, we have
$(\varphi \circ \pi _{\CCI})(z)\leq
(\varphi \circ f_{j}\circ \pi _{\CCI})(z)$
for almost every $z\in \pi _{\CCI}^{-1}f_{j}^{-1}(J(G))$
with respect to
$\tilde{\nu}$.
Since $\tilde{J}(\tilde{f})=
\cup _{j=1}^{m}\Sigma _{m}(j)\cap
\tilde{J}(\tilde{f})$ and
$\Sigma _{m}(j)\cap \tilde{J}(\tilde{f})
\subset \pi _{\CCI}^{-1}(f_{j}^{-1}(J(G)))$
(the latter follows from $\pi _{\CCI}\tilde{f}((w,x))=f_{j}(x)=
f_{j}(\pi _{\CCI}((w,x)))$ for each
$(w,x)\in \Sigma _{m}(j)\cap
\tilde{J}(\tilde{f})$), it follows that
$(\varphi \circ \pi _{\CCI}(z))\leq
(\varphi \circ \pi _{\CCI}\circ \tilde{f})(z)$
for almost every $z\in \tilde{J}(\tilde{f})$ with
respect to $\tilde{\nu}$. Hence, we have shown
Claim 1.

 By Claim 1, we have
$(\varphi \circ \pi _{\CCI}(z))\leq
(\varphi \circ \pi _{\CCI}\circ \tilde{f})(z)$
for almost every $z\in \tilde{J}(\tilde{f})$ with
respect to $\alpha \tilde{\nu}$, where
$\alpha$ is the function in Lemma~\ref{lem2}.
Let $\psi =\varphi \circ \pi _{\CCI}$.
Then, we obtain for each
$n\in \NN $,
$\psi (z)\leq \frac{1}{n}\sum _{j=0}^{n-1}
\psi \circ \tilde{f}^{j}(z)$ for
almost every $z$ with respect to
$\alpha \tilde{\nu}$.
Note that by Lemma~\ref{lem1}-\ref{lem1-3},
the measure $\alpha \tilde{\nu}$ is
$\tilde{f}$-invariant.
Hence, by Birkhoff's ergodic
theorem (see \cite{DGS}),
we have
$\psi (z)\leq \lim _{n\rightarrow \infty}
\frac{1}{n}\sum _{j=0}^{n-1}(\psi \circ
\tilde{f}^{j})(z)$ for almost every
$z$ with respect to
$\alpha \tilde{\nu}$.
Since $\int \psi\alpha \ d\tilde{\nu}
=\int \lim _{n\rightarrow \infty}
\frac{1}{n}\sum _{j=0}^{n-1}(\psi \circ
\tilde{f}^{j})(z)\ \alpha d\tilde{\nu}(z)$,
which follows from Birkhoff's ergodic theorem
again, it follows that
$\psi (z)=\lim _{n\rightarrow \infty}
\frac{1}{n}\sum _{j=0}^{n-1}
(\psi \circ \tilde{f}^{j})(z)$ for
almost every $z$ with respect to
$\alpha \tilde{\nu}$.
Since $\alpha \tilde{\nu}$ is ergodic
(Lemma~\ref{lem1}-\ref{lem1-3}),
then there exists a constant $c$ such that
$\lim _{n\rightarrow \infty}
\frac{1}{n}\sum _{j=0}^{n-1}
(\psi \circ \tilde{f}^{j})(z)=c$ for
almost every $z$ with respect to
$\alpha \tilde{\nu}$.
Hence, it follows that
$\psi (z)=c$ for almost every
$z$ with respect to $\tilde{\nu}$.
Since $\tau$ and $\nu$ are probability
measures, it follows that
$c=1$. Hence,
$\tau =\nu$.
Since $N_{\delta}^{\ast}\nu =\nu$
(Lemma~\ref{convnt}-\ref{convnt1}) and
$\tau$ is $\delta $-conformal (Lemma~\ref{hausconf}),
by Lemma~\ref{cmnt}-\ref{cmnt3} it follows
that
$\nu =\tau$ is a $\delta$-conformal
measure satisfying the
separating condition with respect to
$\{f_{1},\cdots ,f_{m}\}$.
Since supp $\nu =J(G)$ (Main Theorem A),
by Lemma~\ref{nowhere}, it follows that
$f_{i}^{-1}(J(G))\cap f_{j}^{-1}(J(G))$ is
nowhere dense in $f_{j}^{-1}(J(G))$
for each $(i,j)$ with $i\neq j$.

 Hence, we have shown Proposition~\ref{deltasep}.
\end{proof}
\begin{ex}
Let $f_{1}(z)=z^{2}, f_{2}(z)=\frac{z^{2}}{4}$ and
$f_{3}(z)=\frac{z^{2}}{3}$.
Let $G=\langle f_{1},f_{2},f_{3}\rangle$ and
$\tilde{f}:\Sigma _{3}\times \CCI
\rightarrow \Sigma _{3}\times \CCI$ be the
skew product map with respect to
$\{f_{1},f_{2}, f_{3}\}$.
Then, it is easy to see
$J(\langle f_{1},f_{2}\rangle)=
\{z\mid 1\leq |z|\leq 4\}$.
Since $f_{3}^{-1}(J(\langle f_{1},f_{2}\rangle))
=\{z\mid \sqrt{3}\leq |z|\leq 2\sqrt{3}\}
\subset J(\langle f_{1},f_{2}\rangle)$,
we have $J(G)=\{z\mid
1\leq |z|\leq 4.\}$.
Then, $P(G)=\{0,\infty \} \subset F(G)$.
By Theorem 2.6 in \cite{S2}, we find that
$G$ is expanding.
Furthermore, we have $0<H^{2}(J(G))<\infty$ and
$H^{2}(f_{1}^{-1}(J(G))\cap f_{3}^{-1}(J(G)))>0$.
Hence, by Proposition~\ref{deltasep}, 
the number $\delta$ in Lemma~\ref{lem2}
for $\tilde{f}$ satisfies
$\delta >2$.
\end{ex}

\section{Main Theorem B}
\label{secmainB}
In this section, we demonstrate Main Theorem B.
First, we need the following notation.
\begin{df}
Let $G=\langle f_{1},\cdots ,f_{m}\rangle$ be
a finitely generated rational semigroup.
Let $U$ be a non-empty open set in $\CCI$.
We say that $G$ satisfies the {\bf open set condition}
with $U$ with respect to
the generator system $\{f_{1},\cdots ,f_{m}\}$ if
$f_{j}^{-1}(U)\subset U$ for each $j=1,\cdots ,m$ and
$\{f_{j}^{-1}(U)\} _{j=1}^{m}$ are mutually disjoint.
\end{df}

\begin{lem}
\label{oscj}
\begin{enumerate}
\item
If a rational semigroup
$G=\langle f_{1},\cdots ,f_{m}\rangle$
satisfies the open set condition with
$U$ and $\sharp J(G)\geq 3$, then
$J(G)\subset \overline{U}$.
\item If a rational semigroup
$G=\langle f_{1},\cdots ,f_{m}\rangle$
is expanding and if $G$
satisfies the open set condition with
$U$, then $J(G)\subset \overline{U}$.
\end{enumerate}
\end{lem}

\begin{proof}
By Lemma 2.3 (f) in \cite{S5} and Lemma~\ref{elemexp}, 
it is easy to see the statement.
\end{proof}

%

To show Main Theorem B,
we need the following key lemma.
\begin{lem}
\label{upper}
Let
$G=\langle f_{1},\cdots ,f_{m}\rangle$ be a
finitely generated expanding rational semigroup
satisfying the open set condition with
an open set $U$ with respect to $\{f_{1},\cdots ,f_{m}\}$.
Then, we have the following.
\begin{enumerate}
\item There exists a positive constant $C$ such that
for each $r$ with $0<r<$ {\em diam} $\CCI$ and each $x\in J(G)$,
we have $C^{-1}r^{\delta}\leq \nu (B(x,r))\leq Cr^{\delta}$.
Furthermore, $0<H^{\delta}(J(G))<\infty$
and $\dim _{H}(J(G))=\overline{\dim}_{B}(J(G))=\delta$.
\item
Suppose that there exists a $t$-conformal measure $\tau$.
Then, there exists a
positive constant $C_{0}$ such that for any $r$ with
$0<r<$ {\em diam} $\CCI$ and
any $x\in J(G)$, we have
$ C_{0}^{-1}r^{t}\leq \tau (B(x,r))\leq C_{0}r^{t}$.
Furthermore, we have
$0<H^{t}(J(G))<\infty$ and
$\dim _{H}(J(G))=t=\delta$.
Moreover,
$\nu$ and $\tau$ are absolutely continuous
with respect to each other.
\end{enumerate}
\end{lem}
To show this lemma, we need several other 
lemmas(Lemma~\ref{bxr}-Lemma~\ref{preproof}). 
We suppose the assumption of Lemma~\ref{upper}, until
the end of the proof of Lemma~\ref{upper}.

\ 

\noindent {\bf Preparation to show Lemma~\ref{upper}:}
\begin{enumerate}
\item
First, we may assume that $\overline{U}\cap P(G)=\emptyset$.
For, let $V$ be a $\epsilon _{0}$-neighborhood
of $P(G)$ with respect to the hyperbolic metric
on $F(G)$. Then, for each $g\in G$,
we have $g(V)\subset V$, which implies that
$W:=U\setminus \overline{V}$ satisfies
$f_{j}^{-1}(W)\subset W$, for each $j=1,\cdots ,m$,
and $\{f_{j}^{-1}(W)\} _{j}$ are mutually disjoint.
Hence we may assume the above.

 Assuming that $\overline{U}\cap P(G)=\emptyset $,
take a number $\epsilon >0$
such that $B(\overline{U}, 2\epsilon)\cap P(G)=\emptyset$.
Then
for each $y\in \overline{U}$ and any $g\in G$,
we can take well-defined inverse branches of $g^{-1}$
on $B(y,2\epsilon)$.
\item Let $\overline{U}=\sum _{j=1}^{k}K_{j}$
be a measurable partition such that
for each $j=1,\cdots ,k$,
we have int $K_{j}\neq \emptyset$ and
diam $K_{j}\leq \frac{1}{10}\epsilon$.
We take a point $z_{j}\in K_{j}$, for each
$j=1,\cdots ,k$.
\item To show Lemma~\ref{upper}, we may assume that:
for each $j=1,\cdots ,k$ and each $w\in {\cal W}^{\ast}$,
if $\gamma$ is an inverse branch of $f_{w}^{-1}$ on
$B(z_{j},2\epsilon),$\ then
we have
\begin{equation}
\label{imsmall}
\mbox{diam}\gamma (A)\leq (\frac{1}{10})^{|w|}\cdot \mbox{diam}A,
\end{equation}
for each subset $A$ of $B(z_{j},2\epsilon)$.
For, for each $n\in \NN$,
$G_{n}$ (see the notation in section~\ref{rs}) satisfies
$J(G_{n})=J(G)$. Further,
if we use $\nu ^{n}$ to denote the Borel probability measure on
$J(G_{n})=J(G)$ constructed by the
generator system $\{f_{w}\mid |w|=n\}$ of
$G_{n}$, for which the construction method
is the same as that for $\nu$ from
$\{f_{j}\} $, then
$\nu ^{n}$ satisfies $(N_{\delta _{n}}^{n})^{\ast}\nu ^{n}=
\nu ^{n}$ for some $\delta _{n}\in \RR$. Since
$\nu$ satisfies $(N_{\delta}^{n})^{\ast}\nu =\nu $,
by Lemma~\ref{convnt} we obtain $\delta _{n}=\delta$
and $\nu ^{n}=\nu$.
Moreover, since $G$ is expanding,
by the Koebe distortion theorem
there exist numbers $\epsilon '>0$ and
$n\in \NN$ such that
if $\gamma$ is a well-defined inverse
branch of $f_{w}^{-1}$ on $B(z,2\epsilon ')$,
where
$|w|=n$ and $z\in J(G)$, then
for any subset $A$ of $B(z,2\epsilon ')$,
diam $\gamma (A)\leq \frac{1}{10}$ diam $A$.
Let $U':=U\cap B(J(G),\epsilon ')$. Then,
for each $w\in \{1,\cdots ,m\} ^{n}$,
$f_{w}^{-1}(U')\subset U'$ and
$\{f_{w}^{-1}(U')\} _{w: |w|=n}$ are mutually disjoint.
Hence, we may assume the above.
\item Let $r>0$ be fixed. There exists a
number $s\in \NN$ with $s\geq 3$ such that for each
$j=1,\cdots ,k$ and each $w\in {\cal W}^{\ast }$ with 
$|w|\geq s-1$, we have diam $\gamma (K_{j})\leq 
r$,
for each well-defined inverse branch $\gamma$ of
$f_{w}^{-1}$ on $B(z_{j},2\epsilon)$. We fix such an $s$.
Let ${\cal A}$ be the set of all
$(\gamma ,K_{j})$ that satisfies
$j\in \{1,\cdots ,k\}$,
$\gamma$ is a well-defined inverse branch of
$f_{w}^{-1}$ on $B(z_{j},2\epsilon)$ for some
$w\in \{1,\cdots ,m\}^{s}$, and
$\gamma (K_{j})\cap B(x,r)\neq \emptyset$.
\end{enumerate}
Then, we have the following:
\begin{lem}
\label{bxr}
$B(x,r)\cap J(G)=
B(x,r)\cap \bigcup _{(\gamma ,K_{j})\in {\cal A}}
\gamma (J(G)\cap K_{j})$.
\end{lem}
\begin{proof}
Since $J(G)=\cup _{j=1}^{m}f_{j}^{-1}(J(G))$ (Lemma 2.4 in \cite{S5}) and
$J(G)\subset \overline{U}$ (Lemma~\ref{oscj}), it is easy to see the
statement.
\end{proof}

\begin{df}
\begin{enumerate}
\item
Let
$(\gamma ,K_{j})\in {\cal A}$
be any element such that
$\gamma$ is an inverse branch
of $f_{w}^{-1},$\ where
$w=(w_{1},\cdots ,w_{s})\in \{1,\cdots ,m\} ^{s}$.
Then, we decompose $\gamma$
as $\gamma =\gamma _{1}\cdots \gamma _{s}$,
where,
for each $i=1,\cdots ,s$,
we use $\gamma _{i}$ to denote the
inverse branch of $f_{w_{i}}^{-1}$
on $B(\gamma _{i+1}\cdots \gamma _{s}(z_{j}), 2\epsilon)$.
\item For each $A=(\gamma ,K_{j})\in {\cal A}$,
let $l(A)$ be the minimum of $l\in \NN$ that satisfies
$3\leq l\leq s$ and
if $\gamma _{l}\cdots \gamma _{s}(K_{j})\cap K_{i}\neq
\emptyset $, then diam
$\gamma _{1}\cdots \gamma _{l-1}(K_{i})\leq r$.
Note that by (\ref{imsmall}), we have
$\gamma _{1}\cdots \gamma _{l-1}$ is defined on
$K_{i}$ with $\gamma _{l}\cdots \gamma _{s}(K_{j})\cap
K_{i}\neq \emptyset$.
Moreover, note that
according to the choice of $s$,
$l(A)$ exists, for each $A\in {\cal A}$.
\end{enumerate}
\end{df}

\begin{lem}
\label{reqlem}
Let $A=(\gamma ,K_{j})\in {\cal A}$.
If
\begin{equation}
\label{req}
r<\min \{\mbox{{\em diam}}\gamma '_{1}(K_{i})\mid
  (\gamma ',K_{t})\in
{\cal A},\ i\in \{1,\cdots ,k\},\ 
K_{i}\subset B(\gamma '_{2}\cdots \gamma '_{s}(z_{t}),\ 2\epsilon ) \} ,
\end{equation}
then there exists an element
$K_{i}$ such that
$\gamma _{l(A)-1}\cdots \gamma _{s}(K_{j})\cap K_{i}\neq
\emptyset$ and
{\em diam} $\gamma _{1}\cdots \gamma _{l(A)-2}(K_{i})>r$.
\end{lem}

\begin{proof}
If $l(A)\geq 4$, then it is trivial.
If $l(A)=3$, then by (\ref{req}),
the above is true.
\end{proof}

\begin{rem}
\label{remreq}
For the rest, we assume (\ref{req}).
To show Lemma~\ref{upper},
we may make this assumption.
\end{rem}

\begin{df}
For any $A=(\gamma ,K_{j})\in {\cal A}$,
we set
$$\Gamma _{A}:=\{(\gamma _{1}\cdots \gamma _{l(A)-1},\ 
K_{i})\mid \gamma _{l(A)}\cdots \gamma _{s}(K_{j})\cap K_{i}\neq
\emptyset \} .$$
Further, we set $\Gamma =\cup _{A\in {\cal A}}\Gamma _{A}$
(disjoint union).

 Let
$B_{1}$ and $B_{2}$ be two elements in $\Gamma$ with
$B_{1}=(\gamma _{1}\cdots \gamma _{l(A)-1},\ K_{i_{1}})
\in
\Gamma _{A}$ and
$B_{2}=(\gamma '_{1}\cdots \gamma '_{l(A')-1},\ K_{i_{2}})\in
\Gamma _{A'}$, where
$A=(\gamma ,K_{j(A)})\in {\cal A}$ and
$A'=(\gamma ',K_{j(A')})\in {\cal A}$. Then,
\begin{enumerate}
\item We write $B_{1}\sim B_{2}$ if and only if
$K_{i_{1}}=K_{i_{2}}$ and
$\gamma _{1}\cdots \gamma _{l(A)-1}=
\gamma '_{1}\cdots \gamma '_{l(A')-1}$ on
$B(\gamma _{l(A)}\cdots \gamma _{s}
(z_{j(A)}),\ 2\epsilon)\cap
B(\gamma '_{l(A')}\cdots \gamma '_{s}
(z_{j(A')}),\ 2\epsilon)$.
Note that this $\sim$ is an equivalence
relation on $\Gamma $, by
(\ref{imsmall}) and
the uniqueness theorem.
\item We write $B_{1}\preccurlyeq B_{2}$
if and only if
$$\gamma _{1}\cdots \gamma _{l(A)-1}(\mbox{int}
K_{i_{1}})\cap
\gamma '_{1}\cdots \gamma '_{l(A')-1}(\mbox{int}K_{i_{2}})
\neq \emptyset $$
and $l(A)\leq l(A')$.
\end{enumerate}

For any two elements $B$ and $B'$ in
$\Gamma $, we write
$B\preccurlyeq \preccurlyeq B'$ if and only if
there exists a sequence
$\{B_{l}\} _{l=1}^{v}$ in
$\Gamma$ such that
$B=B_{1}\preccurlyeq \cdots \preccurlyeq B_{v}=B'$.
\end{df}

\begin{lem}
\label{rellem1}
Let
$B_{1}$ and $B_{2}$ be two elements in $\Gamma$ with
$B_{1}=(\gamma _{1}\cdots \gamma _{l(A)-1},\ K_{i_{1}})$
$\in
\Gamma _{A}$ and
$B_{2}=(\gamma '_{1}\cdots \gamma '_{l(A')-1},\ K_{i_{2}})\in
\Gamma _{A'}$, where
$A=(\gamma ,K_{j(A)})\in {\cal A}$ and
$A'=(\gamma ',K_{j(A')})\in {\cal A}$.
Suppose that $B_{1}\preccurlyeq B_{2}$.
Then,
we have the following.
\begin{enumerate}
\item \label{rellem1-1}
If $l(A)=l(A')$, then
$B_{1}\sim B_{2}$.
\item \label{rellem1-2}
If $l(A)<l(A')$, then
\begin{enumerate}
\item \label{rellem1-2-1}
{\em int} $K_{i_{1}}\cap \gamma '_{l(A)}\cdots \gamma '_{l(A')-1}($
{\em int} $K_{i_{2}})\neq \emptyset$ and
\item \label{rellem1-2-2}
$\gamma _{1}\cdots \gamma _{l(A)-1}=
\gamma '_{1}\cdots \gamma '_{l(A)-1}$\\ on
$B(\gamma _{l(A)}\cdots \gamma _{s}(z_{j(A)}),\ 2\epsilon)\cap
B(\gamma '_{l(A)}\cdots \gamma '_{s}(z_{j(A')}),\ 2\epsilon)$.
\end{enumerate}
\end{enumerate}
\end{lem}

\begin{proof}
First, we show \ref{rellem1-2}.
Under the assumption of \ref{rellem1-2},
suppose that
$\gamma _{1}\cdots \gamma _{l(A)-1}$ is
an inverse branch of
$f_{w_{1}}^{-1}\cdots f_{w_{l(A)-1}}^{-1}$ and
that $\gamma '_{1}\cdots \gamma '_{l(A')-1}$ is
an inverse branch of
$f_{w'_{1}}^{-1}\cdots f_{w'_{l(A')-1}}^{-1}$.
By the open set condition,
it follows that $w_{j}=w'_{j}$,
for each $j=1,\cdots ,l(A)-1$.
Hence, \ref{rellem1-2-1} holds.

 Next, take a point
$z\in \gamma _{1}\cdots \gamma _{l(A)-1}($ int $K_{i_{1}})\cap
\gamma '_{1}\cdots \gamma '_{l(A')-1}($ int $K_{i_{2}})$.
Let $a:=f_{w_{l(A)-1}}\cdots f_{w_{1}}(z)$.
Then, we have
$a \in$ int $K_{i_{1}}\cap \gamma '_{l(A)}\cdots
\gamma '_{l(A')-1}($ int $K_{i_{2}})$.
Furthermore,
each of $\gamma _{1}\cdots \gamma _{l(A)-1}$ and
$\gamma '_{1}\cdots \gamma '_{l(A)-1}$ is
a well-defined inverse branch of
$(f_{w_{l(A)-1}}\cdots f_{w_{1}})^{-1}$ on
$B(a, \epsilon)$ and
maps $a$ to $z$. Hence, they are equal
on $B(a,\epsilon )$.
By the uniqueness theorem, we obtain
\ref{rellem1-2-2}.

 We can show \ref{rellem1-1} using the same method
as above.
\end{proof}

\begin{lem}
\label{dochi}
If $B$ and $B'$ are two elements
of $\Gamma$ such that
$B\preccurlyeq \preccurlyeq B'$ and
$B'\preccurlyeq \preccurlyeq B$, then
$B\sim B'$.
\end{lem}

\begin{proof}
There exists a sequence
$\{B_{j}\} _{j=1}^{v}$ in $\Gamma$ such that
$B=B_{1}\preccurlyeq \cdots
\preccurlyeq B_{u}=B'\preccurlyeq \cdots
\preccurlyeq B_{v}=B$.
Suppose $B_{j}\in \Gamma _{A_{j}}$, for each
$j=1,\cdots ,v$. Then
we have
$l(A_{1})\leq \cdots \leq l(A_{v})=l(A_{1})$.
By Lemma~\ref{rellem1},
we obtain
$B_{j}\sim B_{j+1}$, for each
$j=1,\cdots ,u-1$.
\end{proof}

\begin{lem}
\label{rellem2}
If $B_{1}\sim B_{2}$,
$B_{3}\sim B_{4}$ and
$B_{1}\preccurlyeq B_{3}$, then
$B_{2}\preccurlyeq B_{4}$.
\end{lem}

\begin{proof}
This is easy to see, from the definitions of
``$\sim $'' and ``$\preccurlyeq $'', by using
(\ref{imsmall}).
\end{proof}
\begin{df}
For any $B\in \Gamma $, we
use $[B]\in \Gamma /\! \sim$ to denote
the equivalence class of $B$,
with respect to the equivalence relation
$\sim$ in $\Gamma$.

 Let $[B_{1}]$ and $[B_{2}]$ be two elements of
$\Gamma /\! \sim $, where
$ B_{1},B_{2}\in \Gamma$.
We write $[B_{1}]\preccurlyeq [B_{2}]$ if and
only if
$B_{1}\preccurlyeq B_{2}$. Note that
this is well defined by
Lemma~\ref{rellem2}. Furthermore,
we write $[B_{1}]\leq [B_{2}]$ if and
only if
$B_{1}\preccurlyeq \preccurlyeq B_{2}$.
Note that this is also well defined
by Lemma~\ref{rellem2}
and that the ``$\leq $'' determines
a partial order in $\Gamma/\! \sim $,
by Lemma~\ref{dochi}.
\end{df}
\begin{lem}
\label{seqb}
Let $q\in \NN$ be an integer with $q\geq 2$.
Let $\{B_{j}\} _{j=1}^{q}$ be a sequence
in $\Gamma$ such that
$B_{1}\preccurlyeq \cdots \preccurlyeq B_{q}$ and
$B_{j}\nsim B_{j+1}$, for each
$j=1,\cdots ,q-1$. Suppose that
for each $j=1,\cdots ,q$,
we have $B_{j}\in \Gamma _{A_{j}}$,
$A_{j}=(\gamma ^{j}, K_{t_{j}})\in {\cal A}$ and
$B_{j}=(\gamma _{1}^{j}\cdots \gamma _{l(A_{j})-1}^{j},\ 
K_{i_{j}})$. Then, we have the following.
\begin{enumerate}
\item \label{seqb1}
$\gamma _{l(A_{1})}^{q}\cdots \gamma _{l(A_{q})-1}^{q}
(K_{i_{q}})\subset B(K_{i_{1}}, \sum\limits _{j=1}^{q-1}
(\frac{1}{10})^{j}\frac{1}{10}\epsilon)$.
\item \label{seqb2}
$\gamma _{1}^{1}\cdots \gamma _{l(A_{1})-1}^{1}=
\gamma _{1}^{q}\cdots \gamma _{l(A_{1})-1}^{q}$\\
on
$V:=B(\gamma _{l(A_{1})}^{1}\cdots
\gamma _{s}^{1}(z_{t_{1}}),\ 2\epsilon)\cap
B(\gamma _{l(A_{1})}^{q}\cdots \gamma _{s}^{q}(z_{t_{q}}),\ 
2\epsilon)$. (Note that by \ref{seqb1}, we have
$V\neq \emptyset$.)
\end{enumerate}
\end{lem}

\begin{proof}
We will show the statement by induction on $q$.
If $q=2$, then the statement follows from
Lemma~\ref{rellem1} and (\ref{imsmall}).
Let $q\geq 3$. Suppose that
the statement holds for each $q'$ with
$2\leq q'\leq q-1$.
By Lemma~\ref{rellem1}, we have
$l(A_{j})<l(A_{j+1})$, for each
$j=1,\cdots ,q-1$.
By the hypothesis of induction,
we have the following claim.\\
Claim 1:
\begin{enumerate}
\item \label{claim1-1}
$\gamma _{l(A_{2})}^{q}\cdots \gamma _{l(A_{q})-1}^{q}
(K_{i_{q}})\subset B(K_{i_{2}}, \sum\limits _{j=1}^{q-2}
(\frac{1}{10})^{j}\frac{1}{10}\epsilon)$.
\item \label{claim1-2}
$\gamma _{1}^{2}\cdots \gamma _{l(A_{2})-1}^{2}=
\gamma _{1}^{q}\cdots \gamma _{l(A_{2})-1}^{q}$\\
on $B(\gamma _{l(A_{2})}^{2}\cdots \gamma _{s}^{2}
(z_{t_{2}}),\ 2\epsilon)\cap
B(\gamma _{l(A_{2})}^{q}\cdots \gamma _{s}^{q}(z_{t_{q}}),\ 
2\epsilon)$.
\end{enumerate}
Combining Claim 1 with (\ref{imsmall}), we obtain
\begin{equation}
\gamma _{l(A_{1})}^{q}\cdots
\gamma _{l(A_{q})-1}^{q}(K_{i_{q}})\subset
B(\gamma _{l(A_{1})}^{2}\cdots
\gamma _{l(A_{2})-1}^{2}(K_{i_{2}}),\ 
\sum\limits _{j=2}^{q-1}(\frac{1}{10})^{j}\frac{1}{10}
\epsilon).
\end{equation}
Moreover, by
Lemma~\ref{rellem1} and (\ref{imsmall}),
we have
$\gamma _{l(A_{1})}^{2}\cdots \gamma _{l(A_{2})-1}^{2}
(K_{i_{2}})\subset B(K_{i_{1}}, \frac{1}{10}\frac{1}{10}
\epsilon)$.
Hence, we obtain
\begin{equation}
\label{gamma1/10}
\gamma _{l(A_{1})}^{q}\cdots \gamma _{l(A_{q})-1}^{q}
(K_{i_{q}})\subset
B(K_{i_{1}}, \sum\limits _{j=1}^{q-1}
(\frac{1}{10})^{j}\frac{1}{10}\epsilon).
\end{equation}
Hence, the statement \ref{seqb1} in our lemma holds for $q$.

 Next, we will show that the statement \ref{seqb2} in our
lemma holds
for $q$.
Let us consider \ref{claim1-2} in Claim 1.
By the open set condition,
for each $j=1,\cdots ,l(A_{1})-1$,
there exists a number $\alpha _{j}\in \{1,\cdots ,m\}$
such that each of $\gamma _{j}^{2}$ and
$\gamma _{j}^{q}$ is an inverse branch of
$f_{\alpha _{j}}^{-1}$. Hence,
we obtain
\begin{equation}
\label{seqbgamma}
\gamma _{l(A_{1})}^{2}\cdots
\gamma _{l(A_{2})-1}^{2}=
\gamma _{l(A_{1})}^{q}\cdots
\gamma _{l(A_{2})-1}^{q}
\end{equation}
on $V_{0}:=
B(\gamma _{l(A_{2})}^{2}\cdots \gamma _{s}^{2}(z_{t_{2}}),\ 
2\epsilon)\cap
B(\gamma _{l(A_{2})}^{q}\cdots \gamma _{s}^{q}(z_{t_{q}}),\ 
2\epsilon)$.

 Let $\beta :=\gamma _{l(A_{1})}^{2}\cdots
\gamma _{l(A_{2})-1}^{2}=
\gamma _{l(A_{1})}^{q}\cdots
\gamma _{l(A_{2})-1}^{q}$ on
$V_{0}$. Then
by \ref{claim1-2} in Claim 1, we obtain
$\gamma _{1}^{2}\cdots \gamma _{l(A_{1})-1}^{2}
=\gamma _{1}^{q}\cdots \gamma _{l(A_{1})-1}^{q}$
on $\beta (V_{0})$. Hence, by the uniqueness theorem,
we get
\begin{equation}
\label{seqb2q}
\gamma _{1}^{2}\cdots
\gamma _{l(A_{1})-1}^{2}=
\gamma _{1}^{q}\cdots
\gamma _{l(A_{1})-1}^{q}
\end{equation}
on $B(\gamma _{l(A_{1})}^{2}\cdots
\gamma _{s}^{2}(z_{t_{2}}),\ 2\epsilon)
\cap
B(\gamma _{l(A_{1})}^{q}\cdots
\gamma _{s}^{q}(z_{t_{q}}),\ 2\epsilon)$.

 Moreover, by Lemma~\ref{rellem1}, we have the following claim.\\
Claim 2:
\begin{enumerate}
\item \label{claim2-1}
int $K_{i_{1}}\cap \gamma _{l(A_{1})}^{2}\cdots
\gamma _{l(A_{2})-1}^{2}($ int $K_{i_{2}})
\neq \emptyset$.
\item \label{claim2-2}
$\gamma _{1}^{1}\cdots
\gamma _{l(A_{1})-1}^{1}=
\gamma _{1}^{2}\cdots
\gamma _{l(A_{1})-1}^{2}$ \\
on
$B(\gamma _{l(A_{1})}^{1}\cdots
\gamma _{s}^{1}(z_{t_{1}}),\ 2\epsilon)
\cap
B(\gamma _{l(A_{1})}^{2}\cdots
\gamma _{s}^{2}(z_{t_{2}}),\ 2\epsilon)$.
\end{enumerate}
Combining \ref{claim2-1} in Claim 2 with
(\ref{imsmall}), we obtain
\begin{equation}
\label{1/5epsilon}
d(\gamma _{l(A_{1})}^{1}\cdots
\gamma _{s}^{1}(z_{t_{1}}),\ 
\gamma _{l(A_{1})}^{2}\cdots
\gamma _{s}^{2}(z_{t_{2}}))
\leq \frac{1}{5}\epsilon.
\end{equation}
Furthermore, by (\ref{gamma1/10})
and (\ref{imsmall}), we obtain
\begin{equation}
\label{3/10epsilon}
d(\gamma _{l(A_{1})}^{q}\cdots
\gamma _{s}^{q}(z_{t_{q}}),\ 
\gamma _{l(A_{1})}^{1}\cdots
\gamma _{s}^{1}(z_{t_{1}}))\leq
\frac{1}{10} \epsilon +
\frac{1}{10} \epsilon +
\sum _{j=1}^{q-1}(\frac{1}{10})^{j}\frac{1}{10} \epsilon
\leq \frac{3}{10} \epsilon.
\end{equation}
Hence, by (\ref{1/5epsilon}) and
(\ref{3/10epsilon}),
we get
$W:=\bigcap _{j=1,2,q}
B(\gamma _{l(A_{1})}^{j}\cdots
\gamma _{s}^{j}(z_{t_{j}}),\ 2\epsilon)\neq \emptyset$.
By \ref{claim2-2} in Claim 2 and (\ref{seqb2q}),
then on $W$,
$\gamma _{1}^{1}\cdots \gamma _{l(A_{1})-1}^{1}
=\gamma _{1}^{q}\cdots \gamma _{l(A_{1})-1}^{q}$.
Hence, by the uniqueness theorem, it follows that
$\gamma _{1}^{1}\cdots
\gamma _{l(A_{1})-1}^{1}=
\gamma _{1}^{q}\cdots
\gamma _{l(A_{1})-1}^{q}$ on
$B(\gamma _{l(A_{1})}^{1}\cdots
\gamma _{s}^{1}(z_{t_{1}}),\ 2\epsilon)
\cap B(\gamma _{l(A_{1})}^{q}\cdots
\gamma _{s}^{q}(z_{t_{q}}),\ 2\epsilon)$.
Hence, the statement \ref{seqb2} in our lemma
holds for
$q$. Hence, the induction is completed.
\end{proof}

\begin{lem}
\label{seqbcor}
Using the same assumption as for
Lemma~\ref{seqb}, it holds that 
$\gamma _{l(A_{1})}^{q}\cdots
\gamma _{s}^{q}(K_{t_{q}})\subset
B(K_{i_{1}}, \frac{1}{5}\epsilon)$.
\end{lem}

\begin{proof}
By Lemma~\ref{seqb} and (\ref{imsmall}), we
obtain
$$
\gamma _{l(A_{1})}^{q}\cdots
\gamma _{s}^{q}(K_{t_{q}})
\subset
B(K_{i_{1}},\ (\sum\limits _{j=1}^{\infty}(\frac{1}{10})^{j})
\frac{1}{10}\epsilon +\frac{1}{10}\epsilon)
\subset
B(K_{i_{1}}, \frac{1}{5}\epsilon).$$
\end{proof}

\begin{df}
Let $\{[m_{1}], \cdots
,[m_{p}]\}$ be the
set of all minimal elements of
$(\Gamma /\! \sim ,\ \leq)$, where,
for each $i=1,\cdots ,p$,
$m_{i}\in \Gamma _{R_{i}},
R_{i}= (\gamma ^{i},K_{u_{i}})\in {\cal A}$ and
$m_{i}=(\gamma _{1}^{i}\cdots \gamma _{l(R_{i})-1}^{i},
K_{v_{i}})$. Furthermore,
for any $i=1,\cdots ,p$,
we use $\eta ^{i}:
\pi _{\CCI}^{-1}(B(\gamma _{l(R_{i})}^{i}\cdots \gamma _{s}^{i}
(z_{u_{i}}),\ 2\epsilon))\rightarrow \Sigma _{m}\times \CCI$ to denote
the inverse branch of $(\tilde{f}^{l(R_{i})-1})^{-1}$
such that $\eta ^{i}((w,y))=
(w^{i}w, \gamma _{1}^{i}\cdots \gamma _{l(R_{i})-1}^{i}(y))$
for each
$(w,y)\in \pi _{\CCI}^{-1}$ $
(B(\gamma _{l(R_{i})}^{i}\cdots \gamma _{s}^{i}
(z_{u_{i}}),\ 2\epsilon))$,
where $w^{i}\in {\cal W}^{\ast}$ is a word
satisfying $|w^{i}|=l(R_{i})-1$
and $\gamma _{1}^{i}\cdots \gamma _{l(R_{i})-1}^{i}$
is an inverse branch of $f_{w^{i}}^{-1}$.
\end{df}

\begin{lem}
\label{preproof}
\begin{enumerate}
\item \label{preproof-1}
$\pi _{\CCI}^{-1}(B(x,r))\cap \tilde{J}(\tilde{f})
\subset
\bigcup _{i=1}^{p}\eta ^{i}(\pi _{\CCI}^{-1}(B(K_{v_{i}},
\frac{1}{5}\epsilon))\cap \tilde{J}(\tilde{f}))$.
\item \label{preproof-2}
$B(x,r)\cap J(G)\subset
\bigcup _{i=1}^{p}\gamma _{1}^{i}\cdots
\gamma _{l(R_{i})-1}^{i}(B(K_{v_{i}}, \frac{1}{5}\epsilon)
\cap J(G))$.
\end{enumerate}
\end{lem}

\begin{proof}
Let $(w,z)\in \pi _{\CCI}^{-1}(B(x,r))\cap
\tilde{J}(\tilde{f})$ be a point.
By Lemma~\ref{oscj}, there exists a number
$j$ such that
$\pi _{\CCI}\tilde{f}^{s}((w,z))\in K_{j}$.
Let $\eta :\pi _{\CCI}^{-1}B(z_{j},2\epsilon)
\rightarrow \Sigma _{m}\times \CCI$ be an
inverse branch of $(\tilde{f}^{s})^{-1}$
such that $\eta ((w',x'))=
((w|s)\cdot w',\ \gamma (x'))$ where
$\gamma$ is an inverse branch of
$f_{w|s}^{-1}$.
Then, we have
$(w,z)\in \eta (\pi _{\CCI}^{-1}(B(z_{j},2\epsilon)))$
and $A:=(\gamma ,K_{j})\in {\cal A}$.
Let $B=(\gamma _{1}\cdots \gamma _{l(A)-1}, K_{i_{1}})
\in \Gamma _{A}$ be an element.
Then, there exists a number $i$ with $1\leq i\leq p$
such that
$[m_{i}]\leq [B]$.
We will show the following claim:\\
Claim 1: $(w,z)\in \eta ^{i}(\pi _{\CCI}^{-1}(B(K_{v_{i}},
\frac{1}{5}\epsilon))\cap \tilde{J}(\tilde{f}))$.

 To show this claim,
we consider the following two cases:\\
Case 1: $B\sim m_{i}$\\
Case 2: There exists a sequence
$(B_{j})_{j=1}^{q}$ in $\Gamma$ such that
$m_{i}=B_{1}\preccurlyeq B_{2}\preccurlyeq
\cdots \preccurlyeq B_{q}=B$ and
$B_{j}\nsim B_{j+1}$ for each
$j=1,\cdots q-1$.

 Suppose that we have Case 2.
Let $y=\pi _{\CCI}(\tilde{f}^{s}((w,z)))\in K_{j}\cap
J(G)$. Then, we have
$z=\gamma (y)=
\gamma _{1}\cdots \gamma _{l(R_{i})-1}\cdot
\gamma _{l(R_{i})}\cdots \gamma _{s}(y)$.
By Lemma~\ref{seqbcor},
we have $\gamma _{l(R_{i})}\cdots
\gamma _{s}(y)\in B(K_{v_{i}},\frac{1}{5}\epsilon)
\cap J(G)$. Furthermore, by
Lemma~\ref{seqb}-\ref{seqb2},
we have
$\gamma _{1}\cdots \gamma _{l(R_{i})-1}=
\gamma _{1}^{i}\cdots \gamma _{l(R_{i})-1}^{i}$
on $B(K_{v_{i}}, \frac{1}{5}\epsilon)$. Combining
this with $B(K_{v_{i}}, \frac{1}{5}\epsilon)\cap
U\neq \emptyset$ and the open set condition,
we get $w|(l(R_{i})-1)=w^{i}$.
By these arguments,
we obtain
$(w,z)= \eta ^{i}(\tilde{f}^{l(R_{i})-1}((w,z)))
\in \eta ^{i}(\pi _{\CCI}^{-1}
(B(K_{v_{i}}, \frac{1}{5}\epsilon)\cap
\tilde{J}(\tilde{f}))$.

 Suppose that we have Case 1. Then,
by the open set condition,
the statement in Claim 1 is true.
Hence, we have shown Claim 1.
 
 By Claim 1, it follows that the statement of
our lemma is true.
\end{proof}
We now demonstrate Lemma~\ref{upper}.\\
{\bf Proof of Lemma~\ref{upper}.}
Let $A=(\gamma ,K_{j})\in {\cal A}$ and
$B=(\gamma _{1}\cdots \gamma _{l(A)-1},\ K_{i})\in
\Gamma _{A}$. By Lemma~\ref{reqlem} and
Remark~\ref{remreq},
there exists a number $u\in \NN$ with
$1\leq u\leq k$ such that
$\gamma _{l(A)-1}\cdots
\gamma _{s}(K_{j})\cap K_{u}\neq \emptyset$ and
diam $\gamma _{1}\cdots \gamma _{l(A)-2}(K_{u})>r$.
Then, by the Koebe distortion theorem,
there exists a positive constant
$C_{1}=C_{1}(\min _{j}\mbox{diam}K_{j}, \epsilon)$,
which is independent of $r,s$ and $x\in J(G)$, such that
$\| (\gamma _{1}\cdots \gamma _{l(A)-2})'(z)\| \geq
C_{1}r$ for each $z\in
B(\gamma _{l(A)-1}\cdots \gamma _{s}(z_{j}),\ \epsilon)$.
Hence, there exists a positive constant
$C_{2}=C_{2}(C_{1},G)$ such that
$$\| (\gamma _{1}\cdots \gamma _{l(A)-1})'(z)\| \geq
C_{2}r,$$ for each
$z\in B(\gamma _{l(A)}\cdots \gamma _{s}(z_{j}),\ \epsilon)$.
Combining this with
$$K_{i}\subset B(\gamma _{l(A)}\cdots \gamma _{s}(z_{j}),\ 
\frac{1}{5}\epsilon ),$$ which
follows from (\ref{imsmall}),
we obtain
$\| (\gamma _{1}\cdots \gamma _{l(A)-1})'(z)||
\geq C_{2}r$, for each
$z\in K_{i}$. Hence, it follows that
there exists a positive constant $C_{3}$, which
is independent of $r,s$ and $x\in J(G)$,
such that
\begin{equation}
\label{intkic3}
\mbox{meas}_{2}(\gamma _{1}\cdots \gamma _{l(A)-1}(
\mbox{int}K_{i}))\geq C_{3}r^{2},
\end{equation}
where meas$_{2}$ denotes the $2$-dimensional
Lebesgue measure.
We now show the following claim:\\
Claim: 
$\gamma _{1}^{i}\cdots \gamma _{l(R_{i})-1}^{i}
(\mbox{int}K_{v_{i}})\subset B(x,3r)$, for each
$i=1,\cdots ,p$.

 To show this claim, since
$\gamma _{1}^{i}\cdots \gamma _{s}^{i}(K_{u_{i}})\cap
B(x,r)\neq \emptyset$ and
$\gamma _{l(R_{i})}^{i}\cdots
\gamma _{s}^{i}(K_{u_{i}})\cap K_{v_{i}}\neq
\emptyset $, we obtain
$\gamma _{1}^{i}\cdots \gamma _{l(R_{i})-1}^{i}(K_{v_{i}})
\cap B(x,2r)\neq \emptyset$. Combining this
with the fact that diam $(\gamma _{1}^{i}\cdots
\gamma _{l(R_{i})-1}^{i}(K_{v_{i}}))\leq r$,
it follows that
the above claim holds.

 Since $\{[m_{1}],\cdots ,[m_{p}]\}$ is the set of
minimal elements of $(\Gamma /\! \! \sim , \leq)$,
we find that
$\{\gamma _{1}^{i}\cdots \gamma _{l(R_{i})-1}^{i}
(\mbox{int}K_{v_{i}})\} _{i=1}^{p}$ are mutually disjoint.
Hence, by (\ref{intkic3}) and the claim, we obtain
\begin{equation}
\label{c4}
p\leq \frac{\mbox{meas}_{2}(B(x,3r))}{C_{3}r^{2}}
\leq C_{4},
\end{equation}
where, $C_{4}$ is a positive constant independent of
$r,s$ and $x\in J(G)$. Furthermore,
by the definition of $l(A)$,
we have diam $\gamma _{1}\cdots \gamma _{l(A)-1}
(K_{i})\leq r$. Hence, by the Koebe distortion theorem,
there exists a positive constant $C_{5}$, which
is independent of $r$ and $x\in J(G)$, such that
\begin{equation}
\label{c5}
\| (\gamma _{1}\cdots \gamma _{l(A)-1})'(z)\|
\leq C_{5}r,
\end{equation}
for each $z\in B(K_{i}, \frac{1}{5}\epsilon)$.
Hence, by Lemma~\ref{preproof},
Lemma~\ref{nuconformal}, Lemma~\ref{cmlemup},
(\ref{c4}) and (\ref{c5}), we obtain
\begin{align*}
\nu (B(x,r))
& = \tilde{\nu}
(\pi _{\CCI}^{-1}(B(x,r))\cap \tilde{J}(\tilde{f}))\\
& \leq \sum\limits _{i=1}^{p}\tilde{\nu}
(\eta ^{i}(\pi _{\CCI}^{-1}
(B(K_{v_{i}}, \frac{1}{5}\epsilon))
\cap \tilde{J}(\tilde{f})))\\
& = \sum\limits _{i=1}^{p}
\int _{\pi _{\CCI}^{-1}(B(K_{v_{i}}, \frac{1}{5}\epsilon))
\cap \tilde{J}(\tilde{f})}
\| (\gamma _{1}^{i}\cdots \gamma _{l(R_{i})-1}^{i})'
(\pi _{\CCI}(z))\| ^{\delta}\ d\tilde{\nu}(z)\\
& \leq C_{4} C_{5}^{\delta}r^{\delta}.
\end{align*}
Similarly, if $\tau$ is a
$t$-conformal measure, then by Lemma~\ref{preproof},
Lemma~\ref{cmlem1}, (\ref{c4}), and (\ref{c5}),
we obtain
\begin{align*}
\tau (B(x,r))
& = \tau (B(x,r)\cap J(G))\\
& \leq \sum\limits _{i=1}^{p}
\tau (\gamma _{1}^{i}\cdots \gamma _{l(R_{i})-1}^{i}
(B(K_{v_{i}}, \frac{1}{5}\epsilon)\cap J(G)))\\
& =\sum\limits _{i=1}^{p}
\int _{B(K_{v_{i}}, \frac{1}{5}\epsilon)\cap J(G)}
\| (\gamma _{1}^{i}\cdots \gamma _{l(R_{i})-1}^{i})'\|
^{t}\ d\tau \\
& \leq C_{4}\cdot C_{5}^{t}\cdot r^{t}.
\end{align*}
By Lemma~\ref{nusubconf} and Lemma~\ref{tsubconf},
we find that
a positive constant
$C'$ exists such that for each
$r$ with $0<r<$ diam $\CCI$ and
$x\in J(G)$,
we have $\nu (B(x,r))\geq C'r^{\delta}$.
Hence, it follows that a
positive constant $C_{6}$ exists such that
for each $r$ with $0<r<$ diam $\CCI$ and
each $x\in J(G)$,
we have $C_{6}^{-1}r^{\delta}
\leq \nu (B(x,r))\leq
C_{6}r^{\delta}$.
Hence, by Proposition 2.2 in \cite{F}
and Main Theorem A,
we obtain
$0<H^{\delta}(J(G))<\infty$ and
$\dim _{H}(J(G))=\overline{\dim}_{B}(J(G))=\delta$.

 Similarly, if $\tau$ is
a $t$-conformal measure, then
by Lemma~\ref{confsubconf},
$\tau$ is $t$-subconformal.
By Lemma~\ref{tsubconf},
we find that
a positive constant $C_{7}$ exists such that
for each
$r>0$ and $x\in J(G)$, we have $\tau (B(x,r))\geq C_{7}r^{t}$.
Hence, it follows that
a positive constant $C_{8}$ exists such that
for each
$r>0$ and $x\in J(G)$, we have
$C_{8}^{-1}r^{t}\leq \tau (B(x,r))\leq C_{8}r^{t}$. Hence,
by Proposition 2.2 in \cite{F},
we obtain
$0<H^{t}(J(G))<\infty$ and
$\dim _{H}(J(G))=t=\delta$.
Then, we find that a positive
constant $C_{9}$ exists such that for each
$x\in J(G)$ and each $r>0$, we have
$C_{9}^{-1}\tau (B(x,r))\leq \nu (B(x,r))
\leq C_{9}\tau (B(x,r))$.
Hence, by the Besicovitch covering lemma
(p294 in \cite{Pe}),
we find that $\nu$ and $\tau$ are absolutely continuous
with respect to each other.
Hence, we have shown Lemma~\ref{upper}.
\qed \\

 We now demonstrate Main Theorem B.\\
{\bf Proof of Main Theorem B:}
By Lemma~\ref{upper}, we find
a positive constant $C$ exists such that
for each $r$ with $0<r<$ diam $\CCI$ and
each $x\in J(G)$, we have
$C^{-1}r^{\delta}\leq \nu (B(x,r))\leq Cr^{\delta}$.
Furthermore,
$\dim _{H}(J(G))=\overline{\dim}_{B}(J(G))=
s(G)=s_{0}(G)=\delta$. Combining this with
Main Theorem A, we see that
for each $x\in \CCI \setminus (A(G)\cup P(G))$, we have
$\dim _{H}(J(G))=S(x)=T(x)=\delta$.

 By Lemma~\ref{upper} and
Proposition~\ref{deltasep},
we obtain
$\nu =\frac{H^{\delta}|_{J(G)}}{H^{\delta}(J(G))}$,
$\nu$ is a $\delta $-conformal measure
satisfying the separating condition
for $\{f_{1},\cdots ,f_{m}\}$, and
$f_{i}^{-1}(J(G))\cap f_{j}^{-1}(J(G))$ is
nowhere dense in $f_{j}^{-1}(J(G))$ for each
$(i,j)$ with $i\neq j$.

 Let $\tau$ be a $t$-conformal measure.
Then, by Lemma~\ref{upper}, we have
$t=\delta$ and $\tau$ is absolutely
continuous with respect to
$\nu$. Since $\nu$ satisfies the
separating condition for $\{f_{1},\cdots
,f_{m}\} $,
it follows that $\tau$ also satisfies the
separating condition for $\{f_{1},\cdots
,f_{m}\}$.
Combining this with Lemma~\ref{uniquecm}-\ref{uniquecm2},
we obtain $\tau =\nu$.

 Hence, we have shown Main Theorem B.
\qed

\section{Examples}
\begin{ex}
\label{hnex1}
\begin{enumerate}
\item
Let $G=\langle f_{1}, f_{2}\rangle$ where $f_{1}(z)=z^{2}$ and $f_{2}(z)
=2.3(z-3)+3$. Then, we can see easily that $\{|z|<0.9\} \subset F(G)$ and
$G$ is expanding.
By the corollary~\ref{cspmaincor}, we get
\[ \overline{\dim}_{B}(J(G))\leq \frac{\log 3}{\log 1.8}<2.\]
In particular, $J(G)$ has no interior points. In \cite{S3}, it was
shown that if a finitely generated rational semigroup satisfies
the open set condition with an open set $U$, then the Julia set is
equal to the closure of the open set $U$ or has no interior points. Note that
the fact that the Julia set of the above semigroup $G$ has no interior
points was shown by using analytic quantity only. It appears to be true
that $G$ does not satisfy the open set condition.
\item Let $G=\langle \frac{z^{3}}{4}, z^{2}+8\rangle$. Then,
we can easily see that $\{|z|<2\} \subset F(G)$ and $G$ is expanding.
Hence, we have
\[ \overline{\dim}_{B}J(G)\leq \frac{\log 5}{\log 3}<2.\]
In particular, $J(G)$ has no interior points.
\end{enumerate}
\end{ex}

\begin{ex}
Let $p_{1},p_{2}$ and $p_{3}\in \CC$ be
mutually distinct points such that
$p_{1}p_{2}p_{3}$ makes a regular triangle.
Let $U$ be the inside part of the regular triangle.
Let $f_{i}(z)=2(z-p_{i})+p_{i}$ for each
$i=1,2,3$.
Let $D(x,r)$ be a Euclidean disk with
radius $r$ in $U\setminus
\cup _{i=1}^{3}f_{i}^{-1}(\overline{U})$,
where $x$ denotes the barycenter of
the regular triangle $p_{1}p_{2}p_{3}$.
Let $g$ be a polynomial such that
$J(g)=\partial D(x,r)$.
Let
$f_{4}(z)=g^{s}(z)$, where $s\in \NN$ is a
large number such that
$f_{4}^{-1}(U)\subset U\setminus \cup _{i=1}^{3}
f_{i}^{-1}(U)$.
Let $G=\langle f_{1},f_{2},f_{3},f_{4}\rangle$.
Then, $G$ satisfies the open set condition
with $U$ with respect to $\{f_{i}\}$.
Furthermore, $G$ is hyperbolic. Hence,
$G$ is expanding, by Theorem 2.6 in \cite{S2}.
Hence, $G$ satisfies the assumption in
Main Theorem B.
(Note that $J(\langle f_{1},f_{2},f_{3}\rangle)$ is
the Sierpi\'{n}ski gasket.)
\end{ex}
\begin{ex}
For any $b$ with $0<b\leq 0.1$,\ 
there exists an $a$ with $0.2<a\leq 1$
 such that $G=\langle 
 a(z-b)^{3}+b,\ z^{2}\rangle $
 satisfies that 
 (1)$G$ is expanding,\ 
 (2)$G$ satisfies the open set condition 
 with $U=\{ z\in \CC \mid |z-b|<\frac{1}{\sqrt{a}},\ |z|>1\} $,\ 
 (3)$J(G)$ is connected,\ and 
 (4)$J(G)$ is porous (hence $\delta =\dim _{H}(J(G))=\overline{\dim }
 _{B}(J(G))<2$). 
\end{ex}

\end{document}